\documentstyle[a4,11pt,bj,etcetera,xy]{article}
\newcommand{\comment}[1]{ }

\newcommand{\iso}{\cong}

\newcommand{\adjshort}[4]{\mbox{$#1,#2:#3\dashv#4$}}

\newcommand{\Set}{\mbox{${\mathcal S}\!et$}}

\newcommand{\morph}[3]{\mbox{$#1:#2\rightarrow #3$}}

\newcommand{\monomorph}[3]{\mbox{$#1:#2\hookrightarrow #3$}}
\newcommand{\cell}[3]{\mbox{$#1:#2\Rightarrow #3$}}

\newcommand{\Subobj}[1]{\mbox{$\mbox{\it Sub}(#1)$}}

\newcommand{\Sub}[1]{\mbox{$\mbox{\sf Sub}(#1)$}}

\newcommand{\twocatdef}[4]{
\setbox153 = \hbox{#1{\bf\hspace*{5mm}morphisms\hspace*{5mm}}}
\catdefsize\textwidth
\advance\catdefsize-\@totalleftmargin
\advance\catdefsize-\wd153
\begin{tabbing}
#1\={\bf\hspace*{5mm}morphisms\hspace*{5mm}}\=\kill
#1\>{\bf\hspace*{5mm}objects}\>\parbox[t]{\catdefsize}{#2} \\
\>{\bf\hspace*{5mm}morphisms}\>\parbox[t]{\catdefsize}{#3} \\
\>{\bf\hspace*{5mm}2-cells}\>\parbox[t]{\catdefsize}{#4}
\end{tabbing}}

\def\sslice{/{\kern-.7ex}/}

\def\Rel{{\mathcal R}\!el}

\newcommand{\alg}[1]{#1
         \mathchoice{\mbox{-}}{\mbox{-}}{\mbox{\scriptsize -}}{\mbox{\tiny -}}
{\mathcal A}lg}

\newcommand{\oftype}{\,:\,}

\def\wrt.{w.r.t.}

\def\ucorner#1{\raisebox{.5ex}{$\ulcorner$}{\kern-.4ex}#1
         {\kern-.4ex}\raisebox{.5ex}{$\urcorner$}}
\def\lcorner#1{\raisebox{-.5ex}{$\llcorner$}{\kern-.4ex}#1
         {\kern-.4ex}\raisebox{-.5ex}{$\lrcorner$}}

\newcommand{\cat}[1]{\mbox{$\Bbb #1$}}  
\newcommand{\Cat}{\mbox{${\mathcal C\/}\!at$}}  \renewcommand{\frac}[2]{ \prooftree #1 \justifies #2
\thickness=0.08em \endprooftree}

\newcommand{\colimit}{\mbox{$\displaystyle\lim_{\rightarrow}$}}
\newcommand{\tuple}[1]
   {\langle #1 \rangle}

\newcommand{\comp}{\,\circ\,}

\def\MonCat{\mbox{${\mathcal M}\!\mathit{on}\Cat$}}

\def\ParCat{\mbox{${\mathcal P}\!\mathit{ar}\Cat$}}
\def\ParMon{\mbox{${\mathcal P}\!\mathit{ar}{\mathcal M}\!\mathit{on}$}}

\newcommand{\Mon}[1]{\mbox{${\mathcal M}\!\mathit{on}$}(#1)}

\newcommand{\ignore}[1]{} 
 \newcommand{\Spn}[1]{{\bold{Spn}}(#1)}

\newcommand{\tensor}{\mathrel{\otimes}}

\newcommand{\Ptl}{\mbox{${\mathsf{Ptl}}$}}

\newcommand{\bold}[1]{\mathbf{#1}} 
\newcommand{\qed}{\hspace*{\fill}$\Box$}
\newcommand{\rel}[3]{{{#1}\colon{#2}\not\rightarrow{#3}}}

\newcommand{\ptlmap}[3]{{{#1}\colon{#2}\rightharpoonup{#3}}}

\input amssym.def
\input{amssym}
\newdimen\proofrulebreadth \proofrulebreadth=.05em
\newdimen\proofdotseparation \proofdotseparation=1.25ex
\newdimen\proofrulebaseline \proofrulebaseline=2ex
\newcount\proofdotnumber \proofdotnumber=3
\let\then\relax
\def\hfi{\hskip0pt plus.0001fil}
\mathchardef\squigto="3A3B
\newif\ifinsideprooftree\insideprooftreefalse
\newif\ifonleftofproofrule\onleftofproofrulefalse
\newif\ifproofdots\proofdotsfalse
\newif\ifdoubleproof\doubleprooffalse
\let\wereinproofbit\relax
\newdimen\shortenproofleft
\newdimen\shortenproofright
\newdimen\proofbelowshift
\newbox\proofabove
\newbox\proofbelow
\newbox\proofrulename
\def\shiftproofbelow{\let\next\relax\afterassignment\setshiftproofbelow\dimen0 }
\def\shiftproofbelowneg{\def\next{\multiply\dimen0 by-1 }%
\afterassignment\setshiftproofbelow\dimen0 }
\def\setshiftproofbelow{\next\proofbelowshift=\dimen0 }
\def\setproofrulebreadth{\proofrulebreadth}

\def\prooftree{
\ifnum	\lastpenalty=1
\then	\unpenalty
\else	\onleftofproofrulefalse
\fi
\ifonleftofproofrule
\else	\ifinsideprooftree
	\then	\hskip.5em plus1fil
	\fi
\fi
\bgroup
\setbox\proofbelow=\hbox{}\setbox\proofrulename=\hbox{}%
\let\justifies\proofover\let\leadsto\proofoverdots\let\Justifies\proofoverdbl
\let\using\proofusing\let\[\prooftree
\ifinsideprooftree\let\]\endprooftree\fi
\proofdotsfalse\doubleprooffalse
\let\thickness\setproofrulebreadth
\let\shiftright\shiftproofbelow \let\shift\shiftproofbelow
\let\shiftleft\shiftproofbelowneg
\let\ifwasinsideprooftree\ifinsideprooftree
\insideprooftreetrue
\setbox\proofabove=\hbox\bgroup$\displaystyle 
\let\wereinproofbit\prooftree
\shortenproofleft=0pt \shortenproofright=0pt \proofbelowshift=0pt
\onleftofproofruletrue\penalty1
}

\def\eproofbit{
\ifx	\wereinproofbit\prooftree
\then	\ifcase	\lastpenalty
	\then	\shortenproofright=0pt	
	\or	\unpenalty\hfil		
	\or	\unpenalty\unskip	
	\else	\shortenproofright=0pt	
	\fi
\fi
\global\dimen0=\shortenproofleft
\global\dimen1=\shortenproofright
\global\dimen2=\proofrulebreadth
\global\dimen3=\proofbelowshift
\global\dimen4=\proofdotseparation
\global\mscount=\proofdotnumber
$\egroup  
\shortenproofleft=\dimen0
\shortenproofright=\dimen1
\proofrulebreadth=\dimen2
\proofbelowshift=\dimen3
\proofdotseparation=\dimen4
\proofdotnumber=\mscount
}

\def\proofover{
\eproofbit 
\setbox\proofbelow=\hbox\bgroup 
\let\wereinproofbit\proofover
$\displaystyle
}%
\def\proofoverdbl{
\eproofbit 
\doubleprooftrue
\setbox\proofbelow=\hbox\bgroup 
\let\wereinproofbit\proofoverdbl
$\displaystyle
}%
\def\proofoverdots{
\eproofbit 
\proofdotstrue
\setbox\proofbelow=\hbox\bgroup 
\let\wereinproofbit\proofoverdots
$\displaystyle
}%
\def\proofusing{
\eproofbit 
\setbox\proofrulename=\hbox\bgroup 
\let\wereinproofbit\proofusing
\kern0.3em$
}

\def\endprooftree{
\eproofbit 
  \dimen5 =0pt
\dimen0=\wd\proofabove \advance\dimen0-\shortenproofleft
\advance\dimen0-\shortenproofright
\dimen1=.5\dimen0 \advance\dimen1-.5\wd\proofbelow
\dimen4=\dimen1
\advance\dimen1\proofbelowshift \advance\dimen4-\proofbelowshift
\ifdim	\dimen1<0pt
\then	\advance\shortenproofleft\dimen1
	\advance\dimen0-\dimen1
	\dimen1=0pt
	\ifdim  \shortenproofleft<0pt
        \then   \setbox\proofabove=\hbox{%
			\kern-\shortenproofleft\unhbox\proofabove}%
                \shortenproofleft=0pt
        \fi
\fi
\ifdim	\dimen4<0pt
\then	\advance\shortenproofright\dimen4
	\advance\dimen0-\dimen4
	\dimen4=0pt
\fi
\ifdim	\shortenproofright<\wd\proofrulename
\then	\shortenproofright=\wd\proofrulename
\fi
\dimen2=\shortenproofleft \advance\dimen2 by\dimen1
\dimen3=\shortenproofright\advance\dimen3 by\dimen4
\ifproofdots
\then
	\dimen6=\shortenproofleft \advance\dimen6 .5\dimen0
	\setbox1=\vbox to\proofdotseparation{\vss\hbox{$\cdot$}\vss}
	\setbox0=\hbox{%
		\kern\dimen6
		$\vcenter to\proofdotnumber\proofdotseparation
			{\leaders\box1\vfill}$%
		\unhbox\proofrulename}%
\else	\dimen6=\fontdimen22\the\textfont2 
	\dimen7=\dimen6
	\advance\dimen6by.5\proofrulebreadth
	\advance\dimen7by-.5\proofrulebreadth
	\setbox0=\hbox{%
		\kern\shortenproofleft
		\ifdoubleproof
		\then	\hbox to\dimen0{%
			$\mathsurround0pt\mathord=\mkern-6mu%
			\cleaders\hbox{$\mkern-2mu=\mkern-2mu$}\hfill
			\mkern-6mu\mathord=$}%
		\else	\vrule height\dimen6 depth-\dimen7 width\dimen0
		\fi
		\unhbox\proofrulename}%
	\ht0=\dimen6 \dp0=-\dimen7
\fi
\let\doll\relax
\ifwasinsideprooftree
\then	\let\VBOX\vbox
\else	\ifmmode\else$\let\doll=$\fi
	\let\VBOX\vcenter
\fi
\VBOX	{\baselineskip\proofrulebaseline \lineskip.2ex
	\expandafter\lineskiplimit\ifproofdots0ex\else-0.6ex\fi
	\hbox	spread\dimen5	{\hfi\unhbox\proofabove\hfi}%
	\hbox{\box0}%
	\hbox	{\kern\dimen2 \box\proofbelow}}\doll%
\global\dimen2=\dimen2
\global\dimen3=\dimen3
\egroup 
\ifonleftofproofrule
\then	\shortenproofleft=\dimen2
\fi
\shortenproofright=\dimen3
\onleftofproofrulefalse
\ifinsideprooftree
\then	\hskip.5em plus 1fil \penalty2
\fi
}

\xyoption{v2}
\xyoption{2cell}
\xyoption{dvips}
\CompileMatrices

\LaTeXdiagrams

\UseAllTwocells

\begin{document}

\begin{center}
{\Large\bf  Paracategories I: Internal Paracategories and Saturated 
Partial Algebras
}                                  
\vskip5mm
\it Claudio Hermida{$^{*}$} $\qquad$ Paulo Mateus{$^{**}$}
          \def\thefootnote{}
          \footnotetext{
CMA, Mathematics Department, IST, 
Lisbon, Portugal.\\
   \hbox{} \hskip-1.8em{*}
e-mail:{\/\tt chermida@math.ist.utl.pt}. The  author acknowledges 
f{i}nancial support from
FCT project Praxis XXI - 18976 - 98.\\
   \hbox{} \hskip-1.8em{**}
e-mail:{\/\tt pmat@math.ist.utl.pt}. 
The author was partially supported by FCT PRAXIS XXI Project 
PRAXIS/P/MAT/10002/1998 Problog, FCT Project FEDER 
POCTI/2001/MAT/37239 Fiblog 
and  FCT grant SFRH/BPD/5625/2001.}
\end{center}

\begin{center}                           
\number\day\space
\ifcase\month\or January\or February\or March\or April\or 
May\or
June\or July\or August\or September\or October\or November\or
December\fi \space\number\year
\end{center}

\begin{abstract}
  Based on the monoid classif{i}er $\Delta$, we give an alternative 
  axiomatization of Freyd's paracategories, which can 
   be interpreted in any  {\em bicategory of partial maps}. 
  Assuming furthermore a free-monoid monad $T$ in our ambient 
category, 
  and coequalisers satisfying some exactness conditions, 
  we give an abstract {\em envelope\/} construction, 
  putting paramonoids (and paracategories) in the more general 
context of \textit{partial algebras\/}. We 
introduce for the latter the crucial notion of
\textit{saturation\/},
which characterises those partial algebras which are isomorphic to
the ones obtained from their enveloping algebras. We also set up a
factorisation system for partial algebras, via epimorphisms and 
(monic) Kleene
morphisms and relate the latter to saturation.
\end{abstract}

\tableofcontents

\comment{
\part{Paracategories I: Internal Paracategories as Saturated Partial 
Algebras}\label{part:paracategories-I}
}

\section{Introduction}
\label{sec:introduction}

There are two primary sources for this work: Freyd's proposal of 
{\em paracategories}, which arise by considering the restricted 
composition structure on an arbitrary subcollection of arrows of a 
category, and  Mateus \etal. 
\cite{MateusSernadasSernadas1999,Mateus2000} notion of {\em 
precategories\/} which provide a setting for probabilistic automata. 
Precategories are graphs endowed with identities and \textit{binary\/}
partial operation of composition, subject to suitable associative and
unit laws. However, all the examples of precategories provided by the
authors (and which we reexamine in \cite{HermidaMateus02b}) are more 
naturally exhibited as paracategories, which provide a richer
framework to study such structures.

The def{i}nition of paracategory in \cite{Freyd96} is phrased in 
elementary terms, using many-sorted f{i}rst-order Horn formulae of 
partial terms, with primitive binary relations $\succeq$
(if the left-hand side is def{i}ned so is the right-hand one and then 
they are equal) and $=$ (Kleene equality, $a = b$ iff 
$a \succeq b$ and $b \succeq a$ \cf. Remark \ref{Kleene-equality}). 
The motivating examples of paracategories in \ibid. are 
bivariant functors and dinatural transformations, and (the one-object 
example) the composition of untyped $\lambda$-terms in normal form.

The one-object case, a so-called {\em paramonoid}, consists of a set 
$M$ and 
$n$-ary partial operations $\ptlmap{\tensor_n}{M^n}{M}$ (for every 
$n$), suitably related, so as to be `partially associative' and 
unitary 
(see  \sectref{sec:ordinary-paracat} for precise details). Hence the 
major differences of a paramonoid with respect to a monoid are: 
\begin{itemize}
\item The {\em partial nature\/} of the operations
\item The weakened associativity axioms. 
\end{itemize}

Categorical algebra obliterates the (somewhat artif{i}cial) 
distinction, traditional in universal algebra, between primitive and 
derived operations. Hence 
for an ordinary monoid, and more generally for a monoid in a monoidal 
category $\cat{C}$,
the whole collection of its operations can be neatly grouped into a 
{\em strong monoidal functor\/} from $\Delta$ into $\cat{C}$, where 
$\Delta$ is the category of f{i}nite ordinals and monotone functions 
(see \sectref{sec:monoids} for details). Our f{i}rst development is 
to 
give an axiomatization of paramonoids in these terms (see 
Def{i}nition 
\ref{def:internal-paramonoid}), and later express a paracategory as a 
particular instance of this abstract notion of paramonoid (see 
Def{i}nition \ref{def:internal-paracategory}).

The next step is the generalisation to this context of Freyd's 
\textbf{{\em envelope\/}} construction, which allows us to consider 
any paramonoid as a subparamonoid of a monoid with a given 
subob{j}ect 
of its carrier. This construction forces us to consider some 
additional structure in our ambient category $\cat{B}$ (which so far 
was only required to have some pullbacks) to carry it out. 
Specif{i}cally, we must assume that $\cat{B}$ admits a free monoid 
monad, as well as coequalisers, with some exactness properties   
(see \sectref{subsec.internal-envelope}). Then, we can construct the 
adjunction between the category of paramonoids and that of 
monoids-with-distinguished-subob{j}ects, Proposition 
\ref{prop:internal-envelope}, which provides in fact a representation 
theorem for the objects of the former.

With the above additional structure on $\cat{B}$, we observe that 
paramonoids are a special case of {\em partial} $T$-algebras. This is
the technical core of the paper (\sectref{sec.lax-algebras}):
we introduce the category of partial algebras for a cartesian
monad (as a special
instance of the notion of lax algebra
\sectref{subsec.partial-alg-monad}), we formulate  the crucial notion
of \textit{saturation\/} for partial algebras 
(\sectref{subsec.saturation})
and we produce the \textit{enveloping algebra\/} construction 
(\sectref{subsec.enveloping-algebra}).
Our main result (Theorem \ref{cor:envelope-partial-algebras})
asserts that the enveloping algebra has the expected universal
property, and characterises saturated partial algebras as those which 
can 
be recovered from their enveloping ones.
We conclude our foray into partial algebras showing that the notion
of Kleene morphism (taken from Freyd's original formulation of
paracategories) makes sense at the partial algebra level, forming
part of a factorisation system in the more general context of
relational algebras
(\sectref{subsec:Kleene-morphisms-factorisation}), and tie this
notion up with saturation (Corollary 
\ref{cor:factorisation-saturation}).

\section{Freyd's paracategories}
\label{sec:ordinary-paracat}

We recall from \cite{Freyd96} the elementary def{i}nitions of 
paramonoid and paracategory and their corresponding morphisms.

\begin{itemize}
    \item  A \textbf{paramonoid} consists of a set $M$ and $n$-ary 
    partial operations $\ptlmap{\tensor_{n}}{M^{n}}{M}$, which we 
    write indistinctly as $[\_]$ for any arity. These operations are 
    subject to the following axioms:
    \begin{enumerate}
        \item  $\ptlmap{\tensor_{0}}{M^{0}(=1)}{M}$ is total
    
        \item  $\ptlmap{\tensor_{1}=(\mathit{id},\mathit{id})}{M}{M}$
    
        \item  \label{sat-paramon} If $[\vec{y}]$ is def{i}ned then 
$[\vec{x}[\vec{y}]\vec{z}] = 
        [\vec{x}\vec{y}\vec{z}]$
    \end{enumerate}
    \noindent where the equality in the last axiom above is Kleene 
    equality (if either side is def{i}ned so is the other and then 
they 
    are equal).

    \item  A \textbf{paracategory} $\cat{C}$ consists of a directed 
    graph    
$C_{0}\stackrel{d}{\leftarrow}C_{1}\stackrel{c}{\rightarrow}C_{0}$ 
    and partial $n$-ary operations 
    $\ptlmap{\comp_{n}}{C_{n}}{C_{1}}$, where
    \[ C_{n} = C_{1}\times_{C_{0}}\cdots\times_{C_{0}} C_{1} = 
    \{(f_{1},\ldots,f_{n})|cf_{i} = df_{i+1}, 1\leq i\leq n-1\} \]
\noindent is the set of composable $n$-tuples of arrows.  They are 
subject to the 
following axioms

\begin{enumerate}
        \item  $\ptlmap{\comp_{0}=\iota}{C_{0}}{C_{1}}$ is total. 
        This yields \textbf{identity arrows} 
$\morph{\mathit{id}_{A}}{A}{A}$ 
        for every object $A\in C_{0}$.
    
        \item  
        $\ptlmap{\comp_{1}=(\mathit{id},\mathit{id})}{C_{1}}{C_{1}}$
    
        \item  If ${\comp}_{n}\vec{y}$ is def{i}ned, then 
    $$ \comp_{m+1+n}(\vec{x},\comp_{k}\vec{y},\vec{z}) =
    \comp_{m+k+n}(\vec{x},\vec{y},\vec{z}) $$
    \noindent where $\vec{x}\in C_{m}$, $\vec{y}\in C_{k}$ and
    $\vec{z}\in C_{n}$.
        
    \end{enumerate}
   
    \item  A \textbf{functor} between paracategories $\cat{C}$ and 
$\cat{D}$ 
    is a morphism of graphs 
    $$(\morph{f_{0}}{C_{0}}{D_{0}},\morph{f_{1}}{C_{1}}{D_{1}})$$ 
    \noindent such that if $[\vec{x}]$ is def{i}ned, then
    $f_{1}[\vec{x}] = [\vec{f_{1}x}]$ (notice that this entails
    preservation of identities).  The functor is called a
    \textbf{Kleene functor} if $[\vec{f_{1}x}] = f_{1}y$  implies
    $[\vec{x}] = y$.
    
    \item A \textbf{subparacategory} of a paracategory $\cat{C}$ is a 
subgraph 
    such that the inclusion is a Kleene functor.
\end{itemize}

We will show in \sectref{subsec:Kleene-morphisms-factorisation}
 that Kleene 
functors form part of a factorisation system, in the more general
context of partial algebras.

A category $\cat{C}$ and a subset of arrows $P\subseteq C_{1}$ 
(including 
the identity arrows) determines a 
subparacategory, to wit, that where $[\vec{x}]$ is def{i}ned if the 
composite  
of the tuple $\vec{x}$ in $\cat{C}$ belongs to $P$. Similarly a 
functor between categories $\morph{F}{\cat{C}}{\cat{D}}$ with 
distinguished subsets of arrows $P\subseteq C_{1}$ and $Q\subseteq 
D_{1}$ such that $f_{1}(P)\subseteq Q$ determines a functor between 
the induced paracategories (see \sectref{sec.lax-algebras} for a 
more general version of this construction in the setting of partial 
algebras 
for a monad).
Let $\ParCat$ denote the category of paracategories and functors 
and $\Cat_{P}$ the category of categories with distinguished subsets 
of arrows (containing the identities) and functors compatible with 
such subsets. The construction above yields a functor 
$\morph{U}{\Cat_{P}}{\ParCat}$.

\begin{proposition}[Enveloping Category \cite{Freyd96}]
    \label{prop:enveloping-category}
\mbox{ }\hfill{ }\\
The functor $\morph{U}{\Cat_{P}}{\ParCat}$ admits a fully 
faithful left adjoint $\mathsf{EC}$.
\end{proposition}
\begin{proof}
    Given a paracategory $\cat{C}$ def{i}ne a category 
    $\mathsf{EC}(\cat{C})$ as follows: let ${\mathsf{F}}(\cat{C})$ be 
the 
free category 
    on the underlying graph of $\cat{C}$ and let 
${\mathsf{EC}}(\cat{C})$ 
    be the quotient of $\mathsf{F}(\cat{C})$  by the relation 
    $$\tuple{x_{1},\ldots,x_{n}} \sim \tuple{[x_{1},\ldots,x_{n}]}$$
    \noindent whenever $[x_{1},\ldots,x_{n}]$ is def{i}ned. Such
    equivalence classes carry a naturally distinguished subset $P$,
    namely $\{\tuple{x}/_{\sim} \,\vert\, x\in\cat{C} \}$.
    This 
    construction extends in an obvious manner to functors of 
    paracategories to yield the desired left adjoint. This left
    adjoint being fully faithful is
equivalent to the unit 
$\morph{\eta}{\cat{C}}{U{\mathsf{EC}}(\cat{C})}$ 
being an isomorphism \cite[\S IV.3,Th.1]{MacLane98}. From 
${\mathsf{EC}}(\cat{C})$ and $P$ we recover a paracategory $\cat{C}'$
whose objects are the equivalence classes of singletons 
$\tuple{x}/_{\sim}$
and whose partial composite is given by
$$[\tuple{x_{1}}/_{\sim},\ldots,\tuple{x_{n}}/_{\sim}] =
\tuple{x_{1},\ldots,x_{n}}/_{\sim}$$

To conclude $\cat{C}\cong\cat{C}'$, we must establish that 
$\tuple{x_{1},\ldots,x_{n}}\bar{\sim}\tuple{y}$ iff $[x_{1},
\ldots,x_{n}] = y$, where $\bar{\sim}$ is the `congruence-closure' of $\sim$. This closure is formed by adding 
reflexivity,  enforcing symmetry  and transitivity, and closing it under 
composition. It is clearly the least-fixed point of a monotone 
operator $F_{\sim}$ on the lattice of relations, $\bar{\sim} = \mu 
F_{\sim}$ and by the finitary nature of this operator, 
the least fixed-point is computed as $\bigcup_{i=0}^{\omega}F^{i}$. 
Thus we can use induction on $i$ to deduce that $$(x_{1},\ldots,x_{n}) 
F^{i} (y) \Longrightarrow [x_{1},\ldots,x_{n}] = y$$
\noindent the only interesting  case being that where the pair is in $F^{i}$ 
by closure under composition: let $\vec{x} =(\vec{l}\vec{x'}\vec{r})F^i y$
by means of

\begin{itemize}
    \item  $\vec{x'} F^j y'$

    \item  $(\vec{l}\  y'\  \vec{r}) F^k y$
\end{itemize}
\noindent for some $j,k < i$. By induction hypothesis, $[\vec{x'}] = 
y'$ and $[\vec{l}\ y'\ \vec{r}] = y$. By the third axiom for
paracategories, 

$$ [\vec{l}\ \vec{x'}\ \vec{r}] = [\vec{l}\ y'\ \vec{r}]= y$$
\qed
\end{proof}

One important consequence of the above enveloping category 
construction is that it explains precisely how every paracategory 
arises, namely by specifying a collection of arrows in a given 
category. 
Any paracategory can be recovered fromits enveloping category.

\begin{examples}
    \label{ex.examples}
    Given the envelope construction, it is straightforward to present
    examples of paracategories as categories with a distinguished
    collection of arrows:
    \begin{enumerate}
        \item  (From \cite{Freyd96}) 
Let $M$ be the monoid of $\beta$/$\eta$-equivalence classes of closed 
untyped lambda terms under the operation
of composition (not application),
given by the combinator traditionally named by $B = \lambda xyz.x(yz)$
(hence the composition of $x$ and $y$ would be $Bxy$).
Consider the subset $N$ of (equivalence classes of) normal terms (that is, those on which the
$\beta$-rule is not applicable). The composition of (representatives 
of) two such is not
necessarily normal(isable), hence we get a paramonoid of equivalence classes 
of normal terms.

        \item  Consider categories $\cat{C}$ and $\cat{D}$ and the
    collection of bivariant functors
    $\morph{T}{\cat{C}^{\mathit{op}}\times\cat{C}}{\cat{D}}$. We 
    obtain a paracategory $\mathsf{DiNat}(\cat{C},\cat{D})$ with 
    such bivariant functors as objects and dinatural
    transformations 
    (\cite[Ch.IX,{\S}4]{MacLane98}) as morphisms. The (pointwise) 
composition of 
    dinatural of transformations is not necessarily dinatural,
    hence the paracategory structure. Note that the ambient
    category is that of bivariant functors whose arrows are
    arbitrary collections of
    morphisms $\{\morph{\alpha_{C}}{T(C,C)}{S(C,C)}\}_{C\in\cat{C}}$
    (no dinaturality required). We explore this example further in 
    \cite{HermidaMateus02b}.
    
        \item  Consider a simply typed $\lambda$-calculus $\cal L$,
    and its collection of $\Set$ valued models $\mathsf{Mod}({\cal
    L})$ (such a model is the same thing as a
    cartesian-closure-preserving-functor from the free ccc
    generated by $\cal L$ into $\Set$ \cite{LambekScott86}). For any 
pair of models
    $M$ and $N$ we have the notion of \textit{logical
    relation\/} between them, which in pedestrian terms is a family of
    relations $\rel{R_{\sigma}}{M(\sigma)}{N(\sigma)}$ indexed by
    the types $\sigma$, satisfying:
    $$ R_{\sigma\times\tau}(x,y) \Longleftrightarrow R_{\sigma}(\pi 
x,\pi y) \wedge R_{\tau}(\pi' x,\pi' y)  $$
    $$ R_{\sigma\Rightarrow\tau}(f,g) \Longleftrightarrow 
    \forall x\oftype M(\sigma), y\oftype N(\sigma).\,\,
    R_{\sigma}(x,y) \Longrightarrow R_{\tau}(fx,gy)$$
    \noindent The usual set-theoretic composition of two such
    logical relations is not necessarily a logical relation. We
    thus obtain a paracategory whose objects are $\Set$-valued
    models of $\cal L$ and whose morphisms are logical relations 
    between such models. As in the previous example, the
    ambient category has models as objects and type-indexed
    collections of relations as arrows (no logical-relation
    condition).
    
\end{enumerate}
\end{examples}

\section{Internal paramonoids and paracategories}
\label{sec:internal-paramonoids}

As we recall in \sectref{sec:monoids}, the monoid classif{i}er 
$\Delta$ (the category of f{i}nite ordinals and monotone maps) gives 
a 
neat way to organise the collection of $n$-ary 
operations of a monoid, so as to guarantee the relevant 
associativity conditions
\[\tensor_{n}(\tensor_{m_{1}},\ldots,\tensor_{m_{n}}) =  
\tensor_{m_{1}+\ldots + m_{n}}  \]
\noindent  This 
equation is forced by functoriality, because in $\Delta$ the ordinal 
$[1]$ 
is a terminal 
object and hence there is only one morphism from $[(m_{1}+\ldots + 
m_{n})]$ into it. 
We take a similar approach to organise 
the $n$-ary operations associated to a paramonoid. The essential 
difference is that composites will not be preserved.  We will 
require instead lax functoriality, so that
\[\tensor_{n}(\tensor_{m_{1}},\ldots,\tensor_{m_{n}}) \leq
\tensor_{m_{1}+\ldots + m_{n}}  \]
\noindent as partial operations from $M^{m_{1}+\ldots + m_{n}}$ to 
$M$. 

The approach outlined above has the immediate advantage of being 
internalisable in any bicategory of partial maps. Thus, we 
set out to internalise paramonoids as a laxif{i}ed version of the 
notion of monoid, that is as `lax internal monoids in a category of 
partial maps'. In addition to laxity, there is a subtlety concerning 
condition
(\ref{sat-paramon}) in the def{i}nition of a paramonoid, which leads 
us 
to introduce below the concept of \textit{saturation\/} for a lax 
functor into the
bicategory of partial maps.

As we recall in \sectref{sec:partial-maps}, our ambient universe for 
paramonoids is a bicategory of partial maps: we consider a
category $\cat{B}$ and a class of monos $\cal M$ in it with the
appropriate closure conditions to make them suitable \textit{domains 
of partial maps\/} and set-up the bicategory $\Ptl_{{\cal 
M}}(\cat{B})$.

Since the hom-categories
in $\Ptl_{{\cal M}}(\cat{B})$ are mere preorders rather than general 
categories, this 
gadget is also referred to as a \textit{locally ordered category}.
Although we would work in this simplif{i}ed context, we
usually keep a more `constructive'  outlook on `proofs of entailment
between predicates' (those specifying the domains of partial maps) in
the back of our minds, so
that the general bicategorical framework is pervasively (albeit
implicitly) present. Of course, the locally ordered situation allows 
us to
simplify the exposition by dispensing with the coherence conditions
pertaining to the structural 2-cells of a lax functor, but it should
be clear that our def{i}nitions are appropriate for the  more general
case as well. 

Let us reexamine condition (\ref{sat-paramon}) of Freyd's 
def{i}nition 
of paramonoid: If $[\vec{y}]\downarrow$  then 
$[\vec{x}[\vec{y}]\vec{z}] = 
    [\vec{x}\vec{y}\vec{z}]$. The Kleene equality amounts to the 
    following  two statements:
    
\begin{enumerate}
    \item  \label{l-to-r} $[\vec{y}]\downarrow$ and
    $[\vec{x}[\vec{y}]\vec{z}]\downarrow$ imply 
    both $[\vec{x}\vec{y}\vec{z}]\downarrow$ and
    $[\vec{x}[\vec{y}]\vec{z}] = 
    [\vec{x}\vec{y}\vec{z}]$.

    \item  \label{r-to-l} $[\vec{y}]\downarrow$ and 
    $[\vec{x}\vec{y}\vec{z}]\downarrow$ imply  
    $[\vec{x}[\vec{y}]\vec{z}]\downarrow$ (and 
$[\vec{x}[\vec{y}]\vec{z}] = 
    [\vec{x}\vec{y}\vec{z}]$)
\end{enumerate}

Bearing in mind that $\tensor_{1} = \mathit{id}$ (by the second axiom 
of paramonoids),
the f{i}rst statement (\ref{l-to-r}) is clearly the laxity condition
\[\tensor_{m+1+n}(\tensor_{1}^{m},\tensor_{k},\tensor_{1}^{n}) \leq
\tensor_{m+k+n}  \]
\noindent where $m = |\vec{x}|$, $k = |\vec{y}|$ and $n = |\vec{z}|$.
The second statement (\ref{r-to-l}) amounts then to the inclusion of 
domains $$ D_{m+k+n} \subseteq
(\tensor_{1}^m\times\tensor_{k}\times\tensor_{1}^n)^{-1}(D_{m+1+n}) 
\cap
(M^m\times{D_{k}}{\times}M^n) $$
\noindent where $D_{n}\subseteq M^n$ is the domain of def{i}nition of
$\tensor_{n}$ and the inverse image is taken along the total maps of 
the operations. This second condition amounts thus to a certain
\textit{saturation\/} of the domains of the partial operations (it
forces more def{i}nedness than that required by laxity alone). Its
abstract counterpart at the level of lax functors into a bicategory
of partial maps is the following\footnote{See \sectref{sec:partial-maps} for notation.}:

\begin{definition}[Saturated lax functor]
    Given a category $\cat{C}$ and a lax functor 
$\morph{F}{\cat{C}}{\Ptl_{{\cal 
M}}(\cat{B})}$ 
consider a pair of
composable morphisms $\morph{f}{X}{Y}$ and $\morph{g}{Y}{Z}$ in 
$\cat{C}$
and the corresponding structural 2-cell $\delta_{f,g}:Fg{\comp}Ff \leq
F(g{\comp}f)$ as displayed
$$
\xymatrixrowsep{2pc}
 \xymatrixcolsep{2pc}
  \let\objectstyle=\scriptstyle
\begin{diagram}
    &&\ar@^{(->}[dddll]_-{d_{gf}}{P_{gf}}\ar@{->}[dddrr]^-{F(gf)}&& \\
&&\ar@^{(->}[dl]^{Ff^{*}d_{g}}\uto^-{^{}\delta^{}}{{Ff}^{*}P_{g}}\drto&&\\
&\ar@^{(->}[dl]^-{d_{f}}{P_{f}}\drto_-{Ff}&&\ar@^{(->}[dl]^{d_{g}}{P_{g}}\drto_-{Fg}& 
\\
    {FX}& &{FY}&& {FZ}
\end{diagram}
$$
\noindent The adjunction\footnote{$\Subobj{X}$ denotes the preorder
of subob{j}ects of $X$.} $\morph{d_{f}{\comp}(\_)\dashv
d_{f}^{*}}{\Subobj{FX}}{\Subobj{P_{f}}}$ induces by transposition a
2-cell (inclusion)
$\morph{\widehat\delta_{f,g}}{{Ff}^{*}P_{g}}{d_{f}^{*}P_{gf}}$.

The lax functor $F$ is called 
 \textbf{saturated} if, for any composable pair $\tuple{f,g}$, the
 2-cell $\widehat\delta_{f,g}$ is an isomorphism.
\end{definition}

More explicitly, the condition of saturation amounts to requiring 
that the
square
$$
\xymatrixrowsep{2pc}
 \xymatrixcolsep{2pc}
\begin{diagram}
{{Ff}^{*}P_{g}}\ar@^{(->}[d]_-{Ff^{*}d_{g}}\ar@^{(->}[r]^-{\delta_{f,g}}&{P_{gf}}
\ar@^{(->}[d]^-{d_{gf}}\\
    {P_{f}}\ar@^{(->}[r]_-{d_{f}}&{FX}
\end{diagram}
$$
\noindent be a pullback, and thus, an intersection of subobjects, as 
alluded in the analysis preceeding the definition.

\begin{remark}\label{rem:saturation-descent}
    We will further examine saturation in the more general context of
    \textit{partial algebras\/} in \sectref{subsec.saturation}, where 
    we unveal its actual meaning as \textit{descent data\/}
    (Proposition \ref{prop:saturation-descent}).
\end{remark}

Recall from \sectref{sec:partial-maps} that $\Ptl_{{\cal M}}(\cat{B})$
has a monoidal structure given by the f{i}nite products of $\cat{B}$.
We are f{i}nally in position to state our main def{i}nition:

\begin{definition}\label{def:internal-paramonoid}
    Consider a category $\cat{B}$ with f{i}nite limits and a class of 
    monos $\cal M$ as in \sectref{sec:partial-maps}
    \begin{itemize}
        \item  An 
    \textbf{internal paramonoid} $\mathbf{M}$ in $\cat{B}$ is a 
\textbf{strong-monoidal 
   saturated lax 
    functor} $\morph{\mathbf{M}}{\Delta}{\Ptl_{{\cal M}}(\cat{B})}$, 
    such that the partial map 
    
$$\ptlmap{F(\morph{\textrm{!`}}{\emptyset}{[1]})}{\mathbf{M}(\emptyset)}{\mathbf{M}(\bold{1})}$$

\noindent is total.

        \item  A \textbf{morphism of internal paramonoids} is a 
monoidal lax 
transformation with \textit{total\/} components
between the corresponding lax functors.
    \end{itemize}

     We have thus the category 
$\ParMon(\cat{B})$ of 
internal paramonoids in $\cat{B}$ (leaving the class $\cal M$ 
implicit).
\end{definition}

 Let us write $M = \mathbf{M}(\bold{1})$ for the underlying object of
 an internal paramonoid. 
 We proceed to unravel the ingredients of the def{i}nition:

\begin{itemize}

    \item  The lax functor $\mathbf{M}$ being strong monoidal 
    implies that it  is determined on objects by $M$ as 
${\mathbf{M}}[n] 
    \iso
    M^{n}$. This is because
    $\Delta$ is 
generated\footnote{The status of $\Delta$ as a monoid classif{i}er
means that it is the \textit{free monoidal category with an internal 
monoid}.} by 
   by $[1]$.

    \item  The (unique) morphism $[n]\rightarrow [1]$ in $\Delta$ 
    yields via $\mathbf{M}$ a partial morphism 
$\ptlmap{\tensor_{n}}{M^{n}}{M}$, whose domain is an $\cal 
M$-subob{j}ect 
$\monomorph{d_n}{D_n}{M^n}$.

    \item  Laxity of $\mathbf{M}$ amounts then to 
    \[\tensor_{n}(\tensor_{m_{1}},\ldots,\tensor_{m_{n}}) \leq
\tensor_{m_{1}+\ldots + m_{n}}  \]
\noindent To see this, notice that the left hand side is the image of 
the composite

\[\xymatrixrowsep{2pc}
 \xymatrixcolsep{2pc}
\begin{diagram}
[(m_{1}+\ldots +
m_{n})-1]\ar@<2ex>[r]_-{\vdots}\ar@<-2ex>[r]&[n-1]\rto&[1]
\end{diagram}
\]
\noindent in $\Delta$, the f{i}rst factor being mapped to 
$\tensor_{m_{1}}\times\ldots\times\tensor_{m_{n}}$ because 
$\mathbf{M}$ is 
strong monoidal.

\item Since $\Ptl_{{\cal M}}(\cat{B})$ is locally ordered and the 
identities (which are total) are maximal,
$\mathit{id} \leq M\mathit{id}$ implies 
$$ M\mathit{id} = \tensor_{1} = \mathit{id}$$
\noindent so that $M$ is {\em normal} (and we get the second condition
of the def{i}nition of paramonoid). The same observation applies to 
lax 
transformations between such functors.

\item The requirement that the partial map
$\morph{e = \mathbf{M}\textrm{!`}}{(\mathbf{M}\emptyset 
\cong)\bold{1}}{M}$ be
total means that $M$ has a constant $e$ which, as we show below, 
behaves as
an always-def{i}ned identity for the partial `multiplications'. 

\item Saturation captures the right-to-left direction of Freyd's
third axiom, namely the inclusion of domains
$$ D_{m+k+n} \subseteq
(\tensor_{1}^m\times\tensor_{k}\times\tensor_{1}^n)^{-1}(D_{m+1+n}) 
\cap
(M^m\times{D_{k}}{\times}M^n) $$
\noindent we mentioned before the definition.
\end{itemize}

\begin{remark}
    It is important to {\em not\/} require that 
$\morph{\mathbf{M}}{\Delta}{\Ptl_{{\cal 
    M}}(\cat{B})}$ preserve the local ordering of 
    $\Delta([m],[n])$. Otherwise, as the second referee pointed out
    to us, composition would be total: $([xy],[z])\leq([x],[yz])$ and
    with $x = e$, the LHS is total and thus $[yz]$ would be always
    defined.
\end{remark}

Notice that for a lax functor 
$\morph{\mathbf{M}}{\Delta}{\Ptl_{{\cal 
M}}(\cat{B})}$, we have a partial 
unit $\ptlmap{e}{1}{M}$ corresponding to the unique morphism 
$\emptyset\rightarrow[1]$ 

\begin{proposition}\label{prop:unit}
    For a strong-monoidal saturated lax functor 
$\morph{\mathbf{M}}{\Delta}{\Ptl_{{\cal 
M}}(\cat{B})}$, if the unique morphism 
$\morph{\textrm{!}}{M}{\bold{1}}$ is a 
strong-epimorphism\footnote{Orthogonal to monos (in ${\cal M}$).
See \cite{Borceux94I} for this property in the context of 
factorisation systems.},
the unit $\ptlmap{e}{1}{M}$ is a total map.
\end{proposition}
\begin{proof}
 Laxity and saturation imply that $\mathit{id} = 
 m\comp(\mathit{id}{\tensor}e)$, where $m$ is the partial binary
 composition. Hence, the right-hand expression is 
total and  
Lemma \ref{lem:cancellation-total-map} implies that 
$\mathit{id}{\tensor}e$ is a total map as well, thus its domain map 
$\mathit{id}\times d$ is an 
isomorphism. The unique diagonal f{i}ll-in in 

\[\xymatrixrowsep{1.5pc} \xymatrixcolsep{2.5pc}
        \begin{diagram} 
{M\times\bold{1}}\dto_{{(\mathit{id}{\times}d)^{-1}}^{}}\rto^-{\pi}&{\bold{1}}\ar@{=}[dd]\ar@{-->}[ddl]\\
{M{\times}D(e)}\dto_-{\pi}& \\
 {D(e)}\ar@^{(->}[r]_-{d}&{\bold{1}}
               \end{diagram}
               \]
shows that the domain of $e$, $\monomorph{d}{D(e)}{\bold{1}}$, is an 
isomorphism as well.

\qed
\end{proof}

Thus we end up with a total unit $e$ 
for the paramonoid as soon as we require that $M$  be
non-empty in an internal sense, \ie. that the unique morphism 
$\morph{\textrm{!}}{M}{\bold{1}}$ 
be a strong-epimorphism.

Next we examine the associativity of the various partial 
operations.

\begin{proposition}
    \label{prop:paramonoid-axioms}
  For a paramonoid $\morph{\mathbf{M}}{\Delta}{\Ptl_{{\cal 
M}}(\cat{B})}$, the following hold:

  \begin{enumerate}
  \item \label{one}
$\tensor_n(x_1,\ldots,x_n) = \tensor_{n+1}(x_1,\ldots,e,\ldots,x_n)$

  \item \label{three}
If $\tensor_k(x_1,.,x_k)$ is def{i}ned then 
$$ \tensor_{m+1+n}(y_1,.,y_m,\tensor_k(x_1,.,x_k),z_1,.,z_n) =
\tensor_{m+k+n}(y_1,.,y_m,x_1,.,x_k,z_1,.,z_n)$$
  \end{enumerate}

\end{proposition}
\begin{proof}
\hfill{ }\\   
\noindent  \underline{(\ref{one})} - 
\begin{eqnarray*}
     \tensor_{n+1}(x_1,\ldots,e,\ldots,x_n)& \succeq & 
\tensor_n(x_1,\ldots,x_n)  \\
     & = & \tensor_n(x_1,\ldots,\tensor(x_i,e),\ldots,x_n)   \\
     & \succeq & \tensor_{n+1}(x_1,\ldots,x_i,e,\ldots,x_n)
\end{eqnarray*}
\noindent where the $\succeq$ follow by laxity, while the equality 
holds by saturation.

\noindent \underline{(\ref{three})} - This is an instance of the
saturation condition:  the left-hand side is
$\tensor_{m+1+k}\comp(\tensor_{1}^{m},\tensor_{k},\tensor_{1}^{n})$.
    Laxity implies the left-to-right inequality. For the
    converse, we must see that the domains of def{i}nition of both
    sides agree. By saturation,
    the intersection of $D_{m+k+n}$ with 
$M^{m}{\times}D_{k}{\times}M^{n}$
    is the same as the preimage of $D_{m+1+n}$ by
    $(\tensor_{1}^{m},\tensor_{k},\tensor_{1}^{n})$, which formulated
    with variables reads as:
    $$ \begin{array}{c}
    (y_1,.,y_m,x_1,.,x_k,z_1,.,z_n) \in D_{m+k+n} \wedge 
(x_1,.,x_k)
    \in D_{k} \\
    \Downarrow \\
    (y_1,.,y_m,\tensor_k(x_1,.,x_k),z_1,.,z_n) \in 
D_{m+1+n}\end{array}$$
\qed
\end{proof}

Now we turn to the consideration of paracategories. Since our 
category $\cat{B}$ has f{i}nite limits, the same is true for 
$\Spn{\cat{B}}(C,C)\equiv \cat{B}/C{\times}C$ and the class of monos 
$\cal M$ can be considered as a suitable class (\textit{dominion\/}) 
for 
$\Spn{\cat{B}}(C,C)$. Hence we can consider the 
bicategory of partial maps $\Ptl_{{\cal M}}(\Spn{\cat{B}}(C,C))$. 
Notice also that a morphism $\morph{f}{C}{D}$ induces, by pullback, a 
{\em change-of-base\/} strong monoidal homomorphism of (monoidal) 
bicategories 
$$\morph{f^{*}}{\Ptl_{{\cal M}}(\Spn{\cat{B}}(D,D))}{\Ptl_{{\cal 
M}}(\Spn{\cat{B}}(C,C))}$$ \noindent which 
takes a paramonoid $\morph{\mathbf{M}}{\Delta}{\Ptl_{{\cal 
M}}(\Spn{\cat{B}}(D,D)))}$ to a paramonoid
$\morph{f^{*}\comp\mathbf{M}}{\Delta}{\Ptl_{{\cal 
M}}(\Spn{\cat{B}}(C,C))}$
(composition with homomorphisms preserves saturation).

\begin{definition}\label{def:internal-paracategory}
  An {\bf internal paracategory\/} in $\cat{B}$ with objects $C_0$ is 
an internal paramonoid $C$
in $\Spn{\cat{B}}(C_0,C_0)$. Given another internal paracategory $D$ 
with objects $D_0$, a \textbf{parafunctor} between them 
is given by a morphism $\morph{f_0}{C_0}{D_0}$ and a morphism of 
internal paramonoids $\morph{f_1}{C}{f_0^{*}(D)}$.
\end{definition}



Let us record the agreement of our def{i}nitions with the elementary 
ones:

\begin{proposition}
    \label{prop:agreement}
  In $\cat{B}= \Set$ with ${\cal M}$ all monos, the def{i}nitions of 
internal paramonoids and paracategories
(Def{i}nitions \ref{def:internal-paramonoid} and 
\ref{def:internal-paracategory}) agree with Freyd's def{i}nitions in
\sectref{sec:ordinary-paracat}, and likewise for their morphisms.
\end{proposition}
\begin{proof}
    Proposition \ref{prop:paramonoid-axioms} shows that our 
def{i}niton
    of internal paramonoid/paracategory specialised to $\Set$
    satisf{i}es Freyd's axioms of \sectref{sec:ordinary-paracat}.
    Conversely, Freyd's third axiom for paramonoid/paracategory
    implies both laxity (left-to-right direction of the Kleene
    equality) and saturation
    (right-to-left).
    \qed
\end{proof}

\begin{remark}
   The above proposition shows that in a sense Freyd's formulation of
   paramonoid is more concise than the one we obtain via lax monoids.
   Yet, our def{i}nition reveals some latent structure in Freyd's
   def{i}niton and puts paramonoids/paracategories in the context of
   \textit{partial algebras\/} 
   (\sectref{sec.lax-algebras}), where we have the right
   abstract ingredients to carry out the envelope construction and 
the representation 
   of saturated partial algebras via their enveloping ones
   (Theorem \ref{cor:envelope-partial-algebras}). 
\end{remark}

\section{The envelope of a paramonoid revisited}
\label{subsec.internal-envelope}

Our next step in the analysis of paramonoids/paracategories is to
give an abstract version of Freyd's enveloping monoid construction,
suitable for our internal treatment. We will reorganise the data of
an internal paramonoid in a more concise form, namely as a
\textit{single\/} partial map, so as to carry out the construction of
the enveloping monoid as a coequaliser. This  links
paramonoids to the partial algebras introduced in 
\sectref{sec.lax-algebras}.

A major advantage of the def{i}nitions of internal paramonoid and 
paracategory 
(Def{i}nitions \ref{def:internal-paramonoid} and 
\ref{def:internal-paracategory} 
respectively) is that they make sense in any category of partial maps 
$\Ptl_{{\cal M}}(\cat{B})$, with a fairly minimal amount of structure 
assumed, 
namely f{i}nite limits in $\cat{B}$ and the closure properties of the 
class $\cal M$ of monos. 

We now need to assume further structure on $\cat{B}$.
Namely, in order to internalise the construction of the enveloping 
monoid/category of Proposition \ref{prop:enveloping-category}, we 
assume that $\cat{B}$ admits the construction of {\em free 
monoids\/}. 
 The usual formula for the free monoid in $\Set$
$\Set$,
\[ TX = \coprod_{n} X^{n} \]
\noindent makes sense in any category $\cat{B}$ with finite products 
such that: 
\begin{enumerate}
    \item  $\cat{B}$ admits countable coproducts.

    \item  The functor $\morph{X\times\_}{\cat{B}}{\cat{B}}$ 
    preserves countable coproducts (for every object $X$).
\end{enumerate}

\begin{remark}
    There are other methods to construct  free monoids. One is by
    taking the colimit of the $\omega$-chain
    \[ 0\rightarrow S0 \rightarrow \cdots S^{n}0 \rightarrow \cdots  
\]
\noindent where  $\morph{S = 
{\mathbf{1}}+(X 
\times \_)}{\cat{B}}{\cat{B}}$. Another alternative is to regard 
$TX$ as an `$X$-labelling' of $T\mathbf{1}$, taking advantage
that the free monoid monad is cartesian, and hence determined by its 
value on the terminal object, $T\mathbf{1}$, which is the
\textit{natural numbers object\/}. The object $T\mathbf{1}$ embodies 
the combinatorics of the free construction; for an arbitrary object 
$X$, $TX$ simply labels the combinatorial data in $T\mathbf{1}$ with 
data from $X$. This works for instance in an
elementary topos with a natural numbers object (which may fail to
have countable coproducts though). See \cite{Benabou90} for a
detailed account.
\end{remark}

In order to construct the envelope of a paramonoid, in addition to
the above conditions which yield free monoids, we assume that 
$\cat{B}$ admits coequalisers stable under
pullbacks (thus $\cat{B}$ is regular \cf. \cite{Borceux94II}).
Let $\Mon{\cat{B}}_{P}$ be the 
category of 
monoids in $\cat{B}$ with a distinguished non-empty subob{j}ect (of 
their 
carriers) and morphisms 
preserving these, and let 
$\morph{U}{\Mon{\cat{B}}_{P}}{\ParMon(\cat{B})}$ be the forgetful 
functor: given a monoid $M$
and a subob{j}ect $\monomorph{d}{D}{M}$, we get partial operations
$\ptlmap{\tensor_{n}}{D^n}{D}$ via the following pullback:

\[\xymatrixrowsep{1.5pc}
 \xymatrixcolsep{2pc}
\begin{diagram}
{\tilde{D}_{n}}\ar@^{(->}[d]_-{d_{n}}
\rrto^-{\tensor_{n}} && {D}\ar@^{(->}[d]^-{d}\\
{D^n}\ar@^{(->}[r]_-{d^n}&{M^n}\ar@^{->}[r]_{\tensor_{n}}&{M}
\end{diagram}
\]
\noindent which organise themselves into a paramonoid structure on
$D$. In $\Set$, $\tilde{D}_{n} =
\{(m_{1},\ldots,m_{n})\in D^n 
\,\vert\,\tensor_{n}(m_{1},\ldots,m_{n})\in D\}$. 
Since the class $\cal M$ of domains of partial maps 
is stable under pullbacks, if $d\in{\cal M}$ then $d^n\in{\cal M}$.
Recall than in any regular category (f{i}nite limits + pullback stable
coequalisers), we can construct the image-factorisation of a morphism
$\morph{f}{X}{Y}$ as follows

\[\xymatrixrowsep{1.5pc} \xymatrixcolsep{2.5pc}
    \begin{diagram}
        &&{\mathit{Im}(f)}\ar@^{(->}[d]^-{m}\\
       {\mathit{Ker}_{f}}
\ar@/^1em/[r]^-{d}\ar@/_1em/[r]_-{c}&{X}\ar@{-->>}[ur]^-{q}\ar@{->}[r]_-{f}&{Y}
\end{diagram}
           \]

\noindent where $\morph{d,c}{\mathit{Ker}_{f}}{X}$ is the
\textit{kernel\/} of $f$, \ie. the pullback of $f$ along itself, and
$\morph{q}{X}{\mathit{Im}(f)}$ is their coequaliser. In $\Set$,
$\mathit{Ker}_{f} = \{(x,y)\,\vert\, fx = fy \}$ and 
${\mathit{Im}}(f) = 
X/\mathit{Ker}_{f}$. The resulting image functor
$\morph{\Sigma_{f}}{\Subobj{X}}{\Subobj{Y}}$ is left-adjoint to the
pullback functor $\morph{f^{*}}{\Subobj{Y}}{\Subobj{X}}$ between the 
categories of subob{j}ects (thus, $\mathit{Im}(f) =
\Sigma_{f}(\mathit{id}_{X})$). Furthermore, such adjoints satisfy a
stability condition, the so-called Beck-Chevalley condition, as a
consequence of the pullback stability of coequalisers. A concise and 
extremely convenient embodiment of these facts is the following:
consider the category $\Sub{\cat{B}}$ whose objects are subob{j}ects
$\monomorph{d}{D}{X}$ and whose morphisms are commutative squares
\[\xymatrixrowsep{1.5pc} \xymatrixcolsep{2.5pc}
    \begin{diagram}
        {D}\ar@^{(->}[d]_-{m}\ar@{-->}[r]&{D'}\ar@^{(->}[d]^-{m'}\\
       {X}\ar@{->}[r]_-{f}&{X'}
           \end{diagram}
\]
          \noindent and let
          $\morph{\mathit{cod}}{\Sub{\cat{B}}}{\cat{B}}$ be the
          evident forgetful functor, taking the codomain of the
          subob{j}ect. The fact that $\cat{B}$ is a regular
          category is equivalent to the statement that $\cat{B}$ 
          has f{i}nite limits and $\mathit{cod}$ is a f{i}bration
          with sums (bif{i}bration satisfying the Beck-Chevalley
          condition) and quotients. See \cite{Jacobs99} for a 
convenient
          account of f{i}bred categorical matters (for the
          logically minded reader).
          
          Since we are working with a distinguished class of
          monos $\cal M$, we would require a similar
          f{i}bration-with-sums structure
          for $\morph{\mathit{cod}}{{\cal M}({\cat{B}})}{\cat{B}}$,
          with ${\cal M}({\cat{B}})$
the evident full subcategory of $\Sub{\cat{B}}$ spanned by the $\cal 
M$-subob{j}ects. We thus have the following requirements of $\cat{B}$
with respect to $\cal M$:

\begin{enumerate}
    \item  $\morph{\mathit{cod}}{{\cal M}({\cat{B}})}{\cat{B}}$ is a 
    f{i}bration with sums.

    \item  The class $\cal M$ is stable under  the relevant colimits 
    (countable coproducts and coequalisers):
    given two  diagrams $\morph{F,G}{\cat{J}}{\cat{B}}$ and a
    natural transformation $\cell{\alpha}{F}{G}$ whose components are
    in $\cal M$, the induced morphism between the colimit objects
    $\morph{{\colimit}\, \alpha}{{\colimit}\, F}{{\colimit}\, G}$ is 
in $\cal M$.
\end{enumerate}

Under these additional assumptions we have:

\begin{proposition}\label{prop:internal-envelope}
    
   The functor $\morph{U}{\Mon{\cat{B}}_{P}}{\ParMon(\cat{B})}$ 
   admits a left adjoint.
\end{proposition}
\begin{proof}
    We outline the construction, bearing in mind that this result is 
    a special case of Theorem \ref{cor:envelope-partial-algebras}.
    Consider a paramonoid $\morph{\mathbf{M}}{\Delta}{\Ptl_{{\cal 
    M}}(\cat{B})}$. 
    We organise all the partial operations 
    $$\ptlmap{\tensor_{n}}{M^{n}}{M} = 
M^{n}\stackrel{d_{n}}{\leftarrow}{D_{n}}\stackrel{s_{n}}{\rightarrow}{M}$$

\noindent into a single span
\[\xymatrixrowsep{1.5pc} \xymatrixcolsep{2.5pc}
        \begin{diagram} 
{TX\iso\coprod_{n}X^{n}}&\lto_-{\coprod_{n}d_{n}}{\coprod_{n}D_{n}}\rto^-{\tuple{s_{n}}}&M  
\end{diagram}
\]
and def{i}ne the {\em enveloping monoid\/} ${{\mathsf{E}}}(M)$ by the 
following coequaliser
\[\xymatrixrowsep{1.5pc} \xymatrixcolsep{3pc}
\let\objectstyle=\scriptstyle
\let\labelstyle=\scriptscriptstyle
        \begin{diagram} 
{T\coprod_{n}D_{n}}\rrcompositemap<-2>_{Td}^{\mu}{\omit} 
\rrto^-{T\tuple{s_{n}}}&&{TM}\ar@{->>}[r]_-{q}&{{\mathsf{E}}(M)}
               \end{diagram}
               \]
         \noindent where $\morph{\mu}{T^{2}M}{TM}$ is the 
multiplication of 
         the free-monoid monad $T$ and $\monomorph{d = 
\coprod_{n}\!d_{n}}{\coprod_{n}D_{n}}{TX}$. Our assumption on the 
stability 
of coequalisers under products(/pullbacks) implies that $T$ preserves 
coequalisers, 
so the upper row of the following diagram is a 
coequaliser:    
\[\xymatrixrowsep{1.5pc} \xymatrixcolsep{2.5pc}
\let\objectstyle=\scriptstyle
\let\labelstyle=\scriptscriptstyle
        \begin{diagram}
            {T^{2}\coprod_{n}\!D_{n}}\dto_-{\mu}
            \ar@/^2ex/[rr]^-{{T^{2}\tuple{s_{n}}}^{}}
\rto_{{T^{2}\coprod_{n}\!d_{n}}^{}}&{T^{3}M}\dto_-{\mu{T}}\rto_-{T\mu}
            &{T^{2}M}\dto^-{\mu}\ar@{->>}[r]^-{Tq}
            &{T{\mathsf{E}}(M)}\ar@{-->}[d]^-{{s_{{\mathsf{E}}(M)}}}\\
 {T\coprod_{n}\!D_{n}}\ar@/_2ex/[rr]_-{{T\tuple{s_{n}}}^{}}
\rto^-{{T\coprod_{n}\!d_{n}}^{}}&{T^{2}M}\rto^-{\mu}&{TM}\ar@{->>}[r]_-{q}&{{\mathsf{E}}(M)}
\end{diagram}
               \]        
\noindent and all squares commute by naturality of $\mu$ and 
associativity for the monad $T$, hence the induced $T$-algebra 
structure on ${\mathsf{E}}(M)$. F{i}nally, ${\mathsf{E}}(M)$ is 
endowed 
with an 
$\cal M$-subob{j}ect, namely  the image  $\Sigma_{q}(d)$.

\qed
\end{proof}

\comment{
If one examines the behaviour of the envelope construction in $\Set$, 
it is not diff{i}cult to observe that it preserves products. In fact, 
we 
can derive this in the more general context set out above under the 
assumption that the monad functor $T$ preserves pullbacks.

\begin{proposition}\label{prop:envelope-preserves-products}
    If $\morph{T}{\cat{B}}{\cat{B}}$ preserves pullbacks, then the 
envelope functor 
    $\morph{{\mathsf{E}}}{\ParMon{\cat{B}}}{\Mon{\cat{B}}_{P}}$ 
    preserves products.
\end{proposition}
\begin{proof}
    Given paramonoids 
$\morph{\mathbf{M}_{1},\mathbf{M}_{2}}{\Delta}{\Ptl_{{\cal 
    M}}(\cat{B})}$, their product is given by their pointwise product 
    in $\cat{B}$. By our hypothesis on $T$ we have a pullback
    
    \[\xymatrixrowsep{1.5pc} \xymatrixcolsep{1.5pc}
        \begin{diagram}
{T(M_{1}\times{M_{2}})}\dto_-{T\pi}\rto^-{{T\pi'}}&{TM_{2}}\dto^-{{T!}}\\
{TM_{1}}\rto_-{T!}&{T\mathbf{1}}
               \end{diagram}
               \]
    \noindent where $T\mathbf{1} = \mathsf{N}$ is the natural numbers 
    object. Since ${T(M_{1}\times{M_{2}})}$ consists of \textit{lists
   of pairs\/},
   the above square identif{i}es 
   \begin{quote}
       \textit{lists of pairs = pairs of
          lists of the same length\/}
   \end{quote}
   \noindent via the canonical comparison 
    morphism 
$\morph{\tuple{T\pi,T\pi'}}{T(M_{1}\times{M_{2}})}{TM_{1}{\times}TM_{2}}$.

Using the units of the monoids, we can transform any pair of lists
   into an equivalent one with both lists of the same length: 
def{i}ne 
    $\morph{\phi}{TM_{1}{\times}TM_{2}}{T(M_{1}\times{M_{2}})}$ by 
    \begin{displaymath}
        \phi(\vec{x},\vec{y}) =  
(\vec{x}.[e_{1}]^{n-\vert\vec{x}\vertÝ},\vec{y}.[e_{2}]^{n-\vert\vec{y}\vert} 

)
    \end{displaymath}
    \noindent where $n= 
    {\mathsf{max}}(\vert\vec{x}\vert ,\vert\vec{y}\vert )$ and 
$e_{1}$ 
    and $e_{2}$ are the units of $M_{1}$ and $M_{2}$
    respectively\footnote{More formally, we could perform this
    construction diagramatically, but the details are unenlightening.
   Alternatively, we could spell out the internal language of our
   ambient category to justify the construction via the above formula
 \cite{Cockett90}  }. 
    Since upon passage to the quotient $(\vec{x},\vec{y}) \sim 
    \phi(\vec{x},\vec{y})$, we conclude 
${{\mathsf{E}}}(M_{1}{\times}M_{2})\iso{{\mathsf{E}}}(M_{1})\times{{\mathsf{E}}}(M_{2})$.

\qed
\end{proof}
}

 %

\section{Partial algebras and saturation}
\label{sec.lax-algebras}

The construction of the envelope of a paramonoid in the proof of 
Proposition \ref{prop:internal-envelope} above prompts us to consider 
the span 
$TM\stackrel{\coprod_{n}d_{n}}{\leftarrow}{T\coprod_{n}D_{n}}\stackrel{\tuple{s_{n}}}{\rightarrow}M$, 

where
$\{\monomorph{d_{n}}{D_{n}}{M}\}_{n}$ is the collection of domains of 
the partial operations of $M$. If the class of monos $\cal M$ is 
closed under coproducts (as it happens for ordinary monos in 
$\Set$), this span is itself a partial map. 
Thus, in the presence of the free-monoid monad $T$, the data for a 
paramonoid can be organised into a {\em single\/} partial map. The 
resulting structure is then an instance of the general notion of 
\textit{partial algebra for a monad\/} which we proceed to examine.

\subsection{Partial algebras for a monad}
\label{subsec.partial-alg-monad}

\begin{definition}
    \label{def:partial-T-algebra}
    Given a category with f{i}nite limits $\cat{B}$ with a class of 
    monos $\cal M$ satisfying the closure conditions of 
    \sectref{sec:partial-maps} and a monad 
    $\tuple{T,\eta,\mu}$ such that $T({\cal M})\subseteq {\cal M}$, a 
    \textbf{partial $T$-algebra} consists of an object $X$ 
    of $\cat{B}$, a partial map $\ptlmap{x}{TX}{X}$ and a 2-cell 
    \[\xymatrixrowsep{1.5pc}
 \xymatrixcolsep{1.5pc}
\begin{diagram}
{T^2X}\ar@^{->}[d]_{{Tx}}
\ar@{}[dr]|-{\stackrel{\alpha}{\leq}}
\rto^{\mu}&{TX}\ar@^{->}[d]^{x}\\
{TX}\ar@^{->}[r]_{x}&{X}
\end{diagram}
\]
satisfying the following \textit{unit condition\/}:
 
$$x\comp\eta = \mathit{id}.$$

Notice that we get also  
$\alpha{\comp}T\eta = \alpha{\comp}\eta_{T} = \mathit{id}$ as there
is only one way in which a partial map ($x$ in this case) is $\leq$ to
itself.

A morphism $\morph{f}{x}{y}$ between partial $T$-algebras 
$\ptlmap{x}{TX}{X}$ and $\ptlmap{y}{TY}{Y}$ is a morphism 
$\morph{f}{X}{Y}$ in $\cat{B}$ together with a 2-cell
   \[\xymatrixrowsep{1.5pc}
 \xymatrixcolsep{1.5pc}
\begin{diagram}
{TX}\ar@^{->}[d]_{{x}}
\ar@{}[dr]|-{\leq}
\rto^-{Tf}&{TY}\ar@^{->}[d]^{y}\\
{X}\ar@{->}[r]_-{f}&{Y}
\end{diagram}
\]
We have thus the category $\Ptl\mbox{\scriptsize -}\alg{T}$ of 
partial $T$-algebras and 
their morphisms.
\end{definition}

A morphism of partial algebras $\morph{f}{x}{y}$ amounts thus to a 
tent
\[\xymatrixrowsep{1.5pc}
 \xymatrixcolsep{1.5pc}
\begin{diagram}
 &\ddlto_{d_{x}}{D_{x}}\drto^{x}\ar@{-->}[rrr]&&& 
\ddlto_(.3){d_{y}}{D_{y}}\drto^{y}&\\
 &&{X}\ar@{->}[rrr]_(.3){f} &&&{Y}\\
 {TX}\ar@{->}[rrr]_-{Tf}&&& {TY}&&   
\end{diagram}
\]

The \textit{unit
condition\/} amounts to the existence of a morphism
$\morph{\eta'_{x}}{X}{D_{x}}$ such that 
$$ d_{x}{\comp}\eta'_{x} = \eta_{X} \qquad x{\comp}\eta'_{x} =
\mathit{id}$$
\noindent and the following being a pullback square

\[\xymatrixrowsep{1.5pc}
 \xymatrixcolsep{1.5pc}
\begin{diagram}
 {X}\dto_-{\mathit{id}}\rto^-{\eta'_{x}}&{D}\dto^-{d_{x}}\\
 {X}\rto_-{\eta_{X}}&{TX}
\end{diagram}
\]
\begin{remark}
    Given $\cat{B}$ with f{i}nite limits and a class of monos $\cal 
M$, 
    consider a monad 
    $\tuple{T,\eta,\mu}$ on ${\cat{B}}$. If the monad is
    \textit{cartesian\/}, we get an 
    induced monad $\morph{\Ptl_{{\cal M}}(T)}{\Ptl_{{\cal 
M}}(\cat{B})}
    {\Ptl_{{\cal M}}(\cat{B})}$ (applying $T$ to a partial map span 
yields 
    another such) with unit $(\mathit{id},\eta)$ and 
    multiplication $(\mathit{id},\mu)$. The {\em lax 
    algebras\/} for $\Ptl_{{\cal M}}(T)$, are the partial algebras
    introduced above. See \cite{Hermida99c} for relevant background
    material on these matters.
\end{remark}

\subsection{Saturation}
\label{subsec.saturation}

In the case of paramonoids, we can recover them from their
enveloping monoids. In the case of partial $T$-algebras, we would
like to recover them from a corresponding \textit{enveloping
algebra\/} construction (see Def{i}nition 
\ref{def:enveloping-algebra} below). It would turn out
that the partial algebras which could be thus obtained would be
characterised (Theorem 
\ref{cor:envelope-partial-algebras}.(\ref{dos}) 
below) by the following crucial condition:

\begin{definition}\label{def:saturated-partial-algebra}
    A partial $T$-algebra 
${TX}\stackrel{d}{\leftarrow}{D}\stackrel{x}{\rightarrow}{X}$ 
    is \textbf{saturated} when the adjoint
    transpose (across $\morph{Td{\comp}(\_)\dashv
(Td)^{*}}{\Subobj{T^{2}X}}{\Subobj{TD}}$) of the structural
    2-cell $\morph{\alpha}{x{\comp}Tx}{x{\comp}\mu}$,
    $\morph{\widehat{\alpha}}{Tx^{*}d}{Td^{*}(\mu^{*}d)}$, is an
    isomorphism. More explicitly, the following square
    $$
    \xymatrixrowsep{2pc}
     \xymatrixcolsep{2pc}
    \begin{diagram}
{{(Tx)}^{*}D}\ar@^{(->}[d]_-{(Tx)^{*}d}\ar@^{(->}[r]^-{\alpha}&{\mu^{*}D}
\ar@^{(->}[d]^-{\mu^{*}d}\\
    {TD}\ar@^{(->}[r]_-{Td}&{T^{2}X}
    \end{diagram}
    $$
 \noindent must be a pullback.   
    
\end{definition}

Recall that in the proof of Proposition \ref{prop:internal-envelope} 
we 
came across the following essential data: the pair of morphisms 
$\morph{Tx,\mu{\comp}Td}{TD}{TX}$. Let us write $\bar{c} = 
\mu{\comp}Td$
and $\bar{d} = Tx$.
There is further data associated 
to the graph
${TX}\stackrel{\bar{d}}{\leftarrow}{TD}\stackrel{\bar{c}}{\rightarrow}{TX}$:

\begin{itemize}
    \item  $\morph{T\eta'_{x}}{TX}{TD}$, which satisf{i}es 
    $$\bar{d}{\comp}T\eta'_{x} = \mathit{id} \qquad 
    \bar{c}{\comp}T\eta'_{x} = \mathit{id}$$

    \item  $\morph{T\alpha}{TD{\bullet}TD}{TD}$, where 
    $$\xymatrixrowsep{2pc}
     \xymatrixcolsep{1.5pc}
    \begin{diagram}
&&\dlto_{\underline{\bar{d}}}{TD{\bullet}TD}\drto^{\underline{\bar{c}}} 
&&\\
    &\dlto_{\bar{d}}TD\drto^{\bar{c}}&
    &\dlto_{\bar{d}}TD\drto^{\bar{c}}& \\
    {TX}& &{TX}& &{TX} 
    \end{diagram}
    $$
    \noindent and the square is a pullback (composition of spans) and
    $T\alpha$ is then a 2-cell between the corresponding spans.
\end{itemize}

\begin{proposition}[Internal category induced by a partial 
$T$-algebra]
    \label{prop:int-cat-ptl-T-alg}
    \hfill{ }\\
    \begin{itemize}
        \item  The data 
$\tuple{(TX,TD,\bar{d},\bar{c}),T\eta'_{x},T\alpha}$
    forms an internal category in $\cat{B}$.
    
        \item  This construction  yields a
    functor $\morph{{\mathsf{cat}}}{\Ptl\mbox{\scriptsize 
-}\alg{T}}{\Cat(\cat{B})}$
    \end{itemize}
    
\end{proposition}

Recall that we have our f{i}bration of $\cal M$-subob{j}ects 
$\morph{\mathit{cod}}{{\cal M}({\cat{B}})}{\cat{B}}$.
We have just seen that a partial algebra $\ptlmap{x}{TX}{X}$ induces 
an internal category ${\mathsf{cat}}(x)$ in $\cat{B}$. Thus it makes 
sense
to consider \textit{descent data\/} for this internal category
in the f{i}bration $\mathit{cod}$. We are not going to make any
substantial use of descent theory as to indulge in details,
and simply refer the reader to
\cite{JanelidzeTholen97,BorceuxJanelidze01}. The point is that the
structural 2-cell $\alpha$ of the partial algebra renders the
following diagram commutative:

\[\xymatrixrowsep{2pc} \xymatrixcolsep{2.5pc}
\let\objectstyle=\scriptstyle
\let\labelstyle=\scriptscriptstyle
    \begin{diagram}
{{(Tx)}^{*}D}\ar@^{(->}[d]_-{(Tx)^{*}d}\ar@/^2ex/[rr]^-{{d^{*}(Tx)}^{}}
\rto_{\alpha^{}}
        &{\mu^{*}D}\ar@^{(->}[d]|-{\mu^{*}d}\rto_-{d^{*}\mu}
        &{D}\ar@^{(->}[d]^-{d} \\
 {TD}\ar@/_2ex/[rr]_-{{Tx}^{}}
 \rto^-{{Td}^{}}
 &{T^{2}X}\rto^-{\mu^{}}&{TX}
\end{diagram}
           \] 
\noindent which is an internal graph in ${\cal M}({\cat{B}})$  (in
fact, an internal discrete f{i}bration) over the category 
${\mathsf{cat}}(x)$. This commutative diagram
induces the morphism 
$\morph{\widehat{\alpha}}{Tx^{*}d}{Td^{*}(\mu^{*}d)}$
which in this context is better read as
$\morph{\widehat{\alpha}}{\bar{d}^{*}d}{\bar{c}^{*}d}$. Now we can
give a different perspective on saturation:

\begin{proposition}\label{prop:saturation-descent}
    The partial $T$-algebra is saturated iff the pair
    
$$(\monomorph{d}{D}{TX},\morph{\widehat{\alpha}}{\bar{c}^{*}d}{\bar{d}^{*}d})$$
    
\noindent constitute descent data for  $\morph{\mathit{cod}}{{\cal 
M}({\cat{B}})}{\cat{B}}$
    over the internal category ${\mathsf{cat}}(x)$.
\end{proposition}
\begin{proof}
    The requirements for descent data are:\hfill{ }
    \begin{itemize}
    \item $\widehat{\alpha}$ must be
    an isomorphism (which is saturation)
    \item  $\widehat{\alpha}$ should satisfy the cocycle conditions, 
    which are automatic in this case, since $d$ (and all its
    pullbacks) is a monomorphism.
    \end{itemize}
    \qed
\end{proof}

As we will see in the proof of Theorem
\ref{cor:envelope-partial-algebras}, this descent-theoretic
perspective
on saturation is helpful in establishing the main characteristic of
saturated partial algebras, namely, that one such is isomorphic to the
partial algebra induced by its  enveloping algebra. As would be clear
in the proof, we could also deduce the relevant properties from
regularity of $\cat{B}$, but the above reformulation of saturation in
terms of descent indicates what we should do if we worked
constructively, keeping track of the `proofs of entailments'.

\subsection{Enveloping Algebra}
\label{subsec.enveloping-algebra}

The construction of the internal enveloping monoid in the 
proof of Proposition \ref{prop:internal-envelope} is already 
formulated at 
the level of partial $T$-algebras: it yields the free $T$-algebra on 
a partial one (plus a distinguished subob{j}ect of its carrier).

\begin{definition}\label{def:T-alg_P}
    Let $\alg{T}_{P}$ 
be the category whose objects $(x,P)$ are  
$T$-algebras $\morph{x}{TX}{X}$ together with a subob{j}ect of the 
carrier 
$\monomorph{m}{P}{X}$
and whose morphisms $\morph{f}{(x,P)}{(y,Q)}$ are 
$T$-algebra morphisms which preserve the subob{j}ects:
\[\xymatrixrowsep{1.5pc}
 \xymatrixcolsep{2pc}
\begin{diagram}
{P}\ar@^{(->}[d]_-{m}
\ar@{-->}[r] & {Q}\ar@^{(->}[d]^-{n}\\
{X}\rto_-{f}&{Y}
\end{diagram}
\]

The 
functor $\morph{U}{\alg{T}_{P}}{\Ptl\mbox{\scriptsize -}\alg{T}}$ 
takes an 
algebra $\morph{x}{TX}{X}$ and a given ($\mathcal{M}$)subob{j}ect 
$\monomorph{m}{D}{X}$ to the partial algebra given by the top span of 
the pullback

 \[\xymatrixrowsep{1.5pc}
 \xymatrixcolsep{2pc}
\begin{diagram}
{\tilde{D}}\ar@^{(->}[d]_-{d}
\rrto^-{x^{\rightharpoonup}} && {D}\ar@^{(->}[d]^-{m}\\
{TD}\ar@^{(->}[r]_-{Tm}&{TX}\ar@{->}[r]_{x}&{X}
\end{diagram}
\]
The axioms for the partial algebra
$\ptlmap{(d,x^{\rightharpoonup})}{TD}{D}$ are easily verif{i}ed: 
\begin{itemize}
    \item  The unit axiom boils down to the fact that the pullback of 
    a mono along itself is the identity span.

    \item  The structural 2-cell
   
$\morph{\alpha}{(d,x^{\rightharpoonup}){\comp}(Td,Tx^{\rightharpoonup})}{(d,x^{\rightharpoonup}){\comp}(\mathit{id},\mu)}$

is induced(via functoriality of pullbacks) by the morphism of spans
$\morph{Tm}{(Tm,Tm)}{(\mathit{id},\mathit{id})}$.

\end{itemize}

\end{definition}

\begin{corollary}[Saturation of induced partial-algebras]
    \label{cor:saturation-induced-partial-algebra}
    \hfill{\mbox{ }}\\
    For any object $(x,P)$, the induced partial algebra $U(x,P)$ is
    saturated.
\end{corollary}
\begin{proof}
    The corresponding structural 2-cell is constructed by pullbacks.
    \qed
\end{proof}

\begin{definition}[Enveloping Algebra]\label{def:enveloping-algebra}
    Assume the ambient category $\cat{B}$ has coequalisers and the 
    monad functor $T$ preserves them. Given a partial $T$-algebra 
$$\ptlmap{x}{TX}{X} = 
    {TX}\stackrel{d_{x}}{\leftarrow}{D}\stackrel{x}{\rightarrow}{X}$$ 
    \noindent def{i}ne its \textbf{enveloping algebra} 
${\mathsf{E}}(x)$
    by the following coequaliser
    \[\xymatrixrowsep{1.5pc} \xymatrixcolsep{3pc}
    \let\objectstyle=\scriptstyle
    \let\labelstyle=\scriptscriptstyle
        \begin{diagram} 
    {TD}\rrcompositemap<-2>_{Td_{x}}^{\mu}{\omit} 
    \rrto^-{Tx}&&{TX}\ar@{->>}[r]^-{q}&{{\mathsf{E}}(x)}
           \end{diagram}
           \]
\noindent with the corresponding algebra structure
$\morph{s_{x}}{T{\mathsf{E}}(x)}{{\mathsf{E}}(x)}$ induced by the
universal property of the coequaliser (preserved by $T$) as in 
the proof of Proposition \ref{prop:internal-envelope}. Likewise, 
we def{i}ne a distinguished subob{j}ect on ${\mathsf{E}}(x)$ by 
taking the
$\cal M$-image of $d_{x}$ along $q$, 
$\monomorph{\Sigma_{q}(d_{x})}{\bar{D}}{{\mathsf{E}}(x)}$.
Thus, the enveloping algebra yields a functor 
$\morph{{\mathsf{E}}}{\Ptl\mbox{\scriptsize -}\alg{T}}{\alg{T}_{P}}$.
\end{definition}

Now we have all the ingredients to state our main result concerning 
partial algebras:

\begin{theorem}\label{cor:envelope-partial-algebras}
    Assume $\cat{B}$ admits, and $T$ preserves, coequalisers, and  
$\cal 
M$
  is closed under images. The following hold:
  \begin{enumerate}
      \item \label{umo}  The functors ${\mathsf E}$ and $U$ are 
adjoint:
      
      \[\xymatrixrowsep{1.5pc} \xymatrixcolsep{3pc}
        \begin{diagram} 
        {\Ptl\mbox{\scriptsize 
-}\alg{T}}\ar@/^2ex/[rr]^-{{\mathsf{E}}}_-{\perp}&&
        \ar@/^2ex/[ll]^-{U}{\alg{T}_{P}}
           \end{diagram}
           \]
\comment{           
      $\morph{{\mathsf{E}}}{\Ptl\mbox{\scriptsize 
-}\alg{T}}{\alg{T}_{P}}$ is left adjoint 
      to
    $\morph{U}{\alg{T}_{P}}{\Ptl\mbox{\scriptsize -}\alg{T}}$.}
  
      \item \label{dos} For a partial algebra $\ptlmap{x}{TX}{X}$, 
the unit of 
      the
    adjunction $\morph{\tilde{\eta}_{x}}{x}{U{\mathsf{E}}(x)}$ is an
    isomorphism iff $\ptlmap{x}{TX}{X}$ is saturated.
  \end{enumerate}
\end{theorem}

\begin{proof}
 \hfill{\mbox{ }}\\   
 \noindent\underline{\ref{umo}}  
To verify that the enveloping algebra has the required universal 
property
$$ \Ptl\mbox{\scriptsize -}\alg{T}(x,U(y,P))\cong
\alg{T}_{P}(({\mathsf{E}}(x),\Sigma_{q}(d_{x})),(y,P))$$

\noindent consider an object $(y,P)$ of $\alg{T}_{P}$, that is, an 
algebra
$\morph{y}{TY}{Y}$ and a subob{j}ect $\monomorph{m}{P}{Y}$, and its
associated partial algebra
\[\xymatrixrowsep{1.5pc}
 \xymatrixcolsep{2pc}
\begin{diagram}
{\tilde{P}}\ar@^{(->}[d]_-{d}
\rrto^-{y^{\rightharpoonup}} && {P}\ar@^{(->}[d]^-{m}\\
{TP}\ar@^{(->}[r]_-{Tm}&{TY}\ar@{->}[r]_-{y}&{Y}
\end{diagram}
\]
and a morphism of partial algebras $\morph{f}{X}{P}$. We get the
required unique morphism $\morph{\hat{f}}{{\mathsf{E}}(x)}{Y}$ as
follows: the morphism $\morph{y{\comp}Tm{\comp}Tf}{TX}{Y}$
coequalises the pair $\morph{(Tx,\mu{\comp}Td_{x})}{TD}{TX}$ and
$\hat{f}$ is the induced factorisation through the coequaliser.

\noindent\underline{\ref{dos}}
We must spell out the unit of the
adjunction $\morph{\tilde{\eta}_{x}}{x}{U{\mathsf{E}}(x)}$:
we have a morphism $\morph{\bar{q}{\comp}\eta'_{x}}{X}{\bar{D}}$,
where $\morph{\bar{q}}{D}{\bar{D}}$ is the instance of the unit of 
the adjunction
$\Sigma_{q} \dashv q^*$ at the object $d$. More explicitly, consider 
the
following diagram:

\[\xymatrixrowsep{1.5pc} \xymatrixcolsep{2.5pc}
\let\objectstyle=\scriptstyle
\let\labelstyle=\scriptscriptstyle
    \begin{diagram}
{{(Tx)}^{*}D}\ar@^{(->}[d]_-{(Tx)^{*}d}\ar@/^2ex/[rr]^-{{d^{*}(Tx)}^{}}
\rto_{\alpha^{}}
    &{\mu^{*}D}\dto_-{\mu^{*}d}\rto_-{d^{*}\mu}
&{D}\ar@^{(->}[d]^-{d}\ar@{->>}[r]^-{\bar{q}}&{\bar{D}}\ar@^{(->}[d]^-{\Sigma_{q}(d)} 
\\
 {TD}\ar@/_2ex/[rr]_-{{Tx}^{}}
 \rto^-{{Td}^{}}
 &{T^{2}X}\rto^-{\mu^{}}&{TX}\ar@{->>}[r]_-{q}&{{\mathsf{E}}(x)}
\end{diagram}
       \] 

\noindent A routine calculation verif{i}es that the morphism 
$\morph{\bar{q}{\comp}\eta'_{x}}{X}{\bar{D}}$
is indeed a morphism of partial $T$-algebras, which is the unit
$\tilde{\eta}_{x}$. 

If $\tilde{\eta}_{x}$ is to be an isomorphism, Corollary 
\ref{cor:saturation-induced-partial-algebra}
says that the partial algebra $x$ must be saturated (saturation is an 
isomorphism-invariant property). It remains
to show that saturation is suff{i}cient. Appealing to
Proposition \ref{prop:saturation-descent}, we descend along $q$ by
making the top row in the above diagram a coequaliser. Thus the image
$\Sigma_{q}(d)$ is computed in this situation as the uniquely induced
morphism between the coequalisers $\bar{q}$ and $q$.  These considerations involving descent could also
be deduced from the axioms of a regular category (specifically, 
coequalisers are pullback stable):
the top part of the above diagram is a pullback of the bottom part
(all squares are pullbacks).

To construct an inverse to $\tilde{\eta}_{x}$, consider the morphism 
$\morph{x}{D}{X}$. Since $\alpha$ is a morphism of partial maps, we
have 

$$ x{\comp}d^*(Tx) = x{\comp}d^*\mu{\comp}\alpha $$

\noindent and by universality of the coequaliser $\bar{q}$, we get a 
unique factorisation of $x$ through $\bar{q}$: there is
$\morph{s}{\bar{D}}{X}$ such that $s{\comp}\bar{q} = x$. This is our 
purported inverse to $\tilde{\eta}_{x}$:

\begin{itemize}
    \item  $s{\comp}\tilde{\eta}_{x} = s{\comp}\bar{q}{\comp}\eta'_{x}
   = x{\comp}\eta'_{x} = \mathit{id} $

    \item  $\tilde{\eta}_{x}{\comp}s = \mathit{id}$ since
    
    \begin{eqnarray*}
        \Sigma_{q}(d){\comp}\tilde{\eta}_{x}{\comp}x & = & 
q{\comp}\eta_{X}{\comp}x  \\
         & = & q{\comp}Tx{\comp}\eta_{D}  \\
         & = & q{\comp}\mu{\comp}Td{\comp}\eta_{D}  \\
         & = & q{\comp}d\\
     & = & \Sigma_{q}(d){\comp}\bar{q}
    \end{eqnarray*}
\noindent and $\Sigma_{q}(d)$ being mono implies 
$$\tilde{\eta}_{x}{\comp}x
= \tilde{\eta}_{x}{\comp}s{\comp}\bar{q} = \bar{q}$$\noindent which
yields the desired equality since $\bar{q}$ is epi.
    \item  We verify that $\morph{s}{\bar{D}}{X}$ is a morphism of
    partial algebras using the abovementioned fact that the square
    with the $q$'s is a pullback and a routine calculation 
    like the previous ones.
\end{itemize}
\qed
\end{proof}

The second assertion in the above theorem means that the saturation 
condition  distinguishes those partial algebras which can be
recovered from their enveloping algebra. In particular, when
$\mathsf{M}$ is
the free-monoid monad on $\Set$ and $\cal M$ is the class of all
(regular) monos, we have our desired identif{i}cation:

\begin{center}
\fbox{$ {\mathsf{Sat^{\bullet}}}(\Ptl-\alg{{\mathsf{M}}}) \simeq 
\ParMon $}
\end{center}

\noindent where ${\mathsf{Sat^{\bullet}}}(\Ptl-\alg{\mathsf{M}})$ is 
the full subcategory 
of saturated partial $\mathsf{M}$-algebras with non-empty carrier $X$ 
($\morph{\textrm{!}}{X}{\bold{1}}$ being a strong-epimorphism). Thus 
the only additional
ingredient present in the case of paramonoids/paracategories is the
assumption that the carrier is non-empty. Clearly if we impose a
similar non-emptiness condition on subob{j}ects in $\alg{T}_{P}$, 
Theorem 
\ref{cor:envelope-partial-algebras} still applies , all the
relevant constructions remaining unchanged.

\begin{corollary}[Reflectivity of saturated partial algebras]
    \label{cor:reflectivity-saturated-ptl-alg}
    \hfill{\mbox{ }}\\
    The full subcategory ${\mathsf{Sat}}(\Ptl\mbox{\scriptsize 
-}\alg{T})$ of  
    saturated partial $T$-algebras is reflective, \ie. the inclusion 
    $\monomorph{\iota}{{\mathsf{Sat}}(\Ptl\mbox{\scriptsize 
-}\alg{T})}{\Ptl\mbox{\scriptsize -}\alg{T}}$
    has a left adjoint \linebreak
    $\morph{{\mathsf{Sat}}}{\Ptl\mbox{\scriptsize 
-}\alg{T}}{{\mathsf{Sat}}(\Ptl\mbox{\scriptsize -}\alg{T})}$.
\end{corollary}
\begin{proof}
    Given an adjunction $\adjshort{F}{G}{\cat{D}}{\cat{C}}$ with unit
   $\eta$ such that ${\eta}G$ is an isomorphism , the full 
    subcategory $\cat{C}'$ of $\cat{C}$ of those objects for which 
the 
    unit is an isomorphism is reflective. Therefore the result 
follows 
    from the caractherisation of saturated partial 
    algebras as those for which the unit of the adjunction of Theorem
    \ref{cor:envelope-partial-algebras} is an isomorphism. In view of 
the
    developments to follow 
    it is instructive to spell out the construction, assuming the
    unit $\cell{\eta}{\mathit{id}}{T}$ is cartesian (see  
    \sectref{subsec:Kleene-morphisms-factorisation}).
    Set $S(x) = U{\mathsf{E}}(x)$ and pull it back along
    $\tilde{\eta}_{x}$: if
    \[\xymatrixrowsep{1.5pc}
     \xymatrixcolsep{1.5pc}
    \begin{diagram}
&\ar@{-->}[ddl]_-{d'}{D_{{\mathsf{Sat}}(x)}}\ar@{-->}[dr]^{x'}\ar@{-->}[rrr]&&& 
\ddlto_(.3){d_{S(x)}}{D_{S(x)}}\drto&\\
     &&{X}\ar@{->}[rrr]_(.3){\tilde{\eta}_{x}} 
    &&&{\bar{D}}\\
     {TX}\ar@{->}[rrr]_-{T\tilde{\eta}_{x}}&&& {T\bar{D}}&&   
    \end{diagram}
    \]
    \noindent is a limit diagram, then ${\mathsf{Sat}}(x) =
    \ptlmap{(d',x')}{TX}{X}$ is the required reflection.
    \qed
\end{proof}

\subsection{Kleene morphisms and a factorisation system for partial 
algebras}
\label{subsec:Kleene-morphisms-factorisation}

In the context of paracategories, Kleene functors play a role in
def{i}ning the relevant notion of subparacategory 
(\sectref{sec:ordinary-paracat}).
Since they are only used in that context, it makes sense to consider 
the `Kleene functor' condition only in the case when the underlying
graph morphism is monic. We will show that such monic Kleene functors
make sense at the level of partial algebras, where they form part of 
a factorisation system. We will also make clear at this level of
generality Freyd's statement that the envelope construction embeds
every paracategory as a subparacategory of a category. This will give
yet another condition equivalent to saturation in Corollary 
\ref{cor:factorisation-saturation}

The notion of monic Kleene morphism can be formulated at the level of 
partial algebras:

\begin{definition}\label{def:Kleene-morphism}
    An $\mathcal{M}$-monomorphism $\monomorph{f}{x}{y}$ between 
partial $T$-algebras 
$\ptlmap{x}{TX}{X}$ and $\ptlmap{y}{TY}{Y}$ is a \textbf{Kleene
morphism} if the corresponding tent diagram
\[\xymatrixrowsep{1.5pc}
 \xymatrixcolsep{1.5pc}
\begin{diagram}
 &\ar@{-->}[ddl]_-{d_{x}}{D_{x}}\ar@{-->}[dr]^-{x}\ar@{-->}[rrr]&&& 
\ddlto_(.3){d_{y}}{D_{y}}\drto^{y}&\\
 &&{X}\ar@{->}[rrr]_(.3){f} 
&&&{Y}\\
 {TX}\ar@{->}[rrr]_-{Tf}&&& {TY}&&   
\end{diagram}
\]
\noindent is a limit diagram, \ie. the dashed arrows constitute a
(co)universal cone for the diagram given by the straight arrows.
\end{definition}

In $\Set$-theoretic formulation: 
$$ \tensor^y(\vec{fy}) = fz \Longrightarrow (\vec{y}\in D_{x} \wedge 
\tensor^x(\vec{y}) = z)$$

\begin{proposition}\label{prop:Kleene-morphisms}
    A functor $\morph{(f_{0},f_{1})}{\cat{C}}{\cat{D}}$ between
    paracategories with $f_{0}\in\mathcal{M}$ is a Kleene functor iff 
the corresponding morphism
    of paramonoids  $\morph{f_1}{C}{f_0^{*}(D)}$ is a Kleene morphism.
\end{proposition}

\begin{proposition}[Fibrational property  of partial algebras]
\label{prop:fib-ptl-alg}
\hfill{\mbox{  }}\\
The forgetful functor 
$\morph{\underline{U}}{\Rel\mbox{\scriptsize -}\alg{T}}{\cat{B}}$, 
which takes 
a partial algebra $\ptlmap{x}{TX}{X}$ to $X$, admits cartesian 
liftings 
of $\mathcal{M}$-monomorphisms. An $\mathcal{M}$-monomorphism of 
partial algebras is cartesian iff it is a Kleene morphism.
Furthermore, if the codomain of a cartesian morphism is saturated so 
is its domain.
\end{proposition}
\begin{proof}
Given an $\mathcal{M}$-monomorphism $\monomorph{f}{X}{Y}$, we can 
construct the cartesian lifting of 
at a partial algebra $\ptlmap{y}{TY}{Y}$ via the above 
 limit 
diagram, which can be computed via pullbacks:


 \[\xymatrixrowsep{1.5pc}
 \xymatrixcolsep{1.5pc}
  \let\objectstyle=\scriptstyle
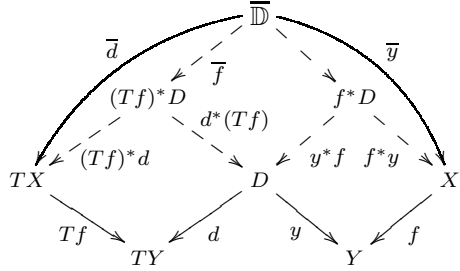
\begin{diagram}  
  &&  
\ar@/_3ex/[ddll]_-{\overline{d}}\ar@{-->}[dl]^-{\overline{f}}{\overline{\cat{D}}}\ar@{-->}[dr]\ar@/^3ex/[ddrr]^-{\overline{y}} 
& & \\
&\ar@{-->}[dl]^-{(Tf)^{*}d}{(Tf)^{*}D}\ar@{-->}[dr]^{d^{*}(Tf)} &  
&\ar@{-->}[dl]^-{y^{*}f}{f^{*}D}\ar@{-->}[dr]_-{f^{*}y} & \\    
{TX}\drto_{Tf}& &\dlto^{d}{D}\drto_{y} & &\dlto^{f}{X}\\    
    &{TY} & &{Y} &   
\end{diagram}
\]

\noindent Notice by the stability properties of $\mathcal{M}$ both 
$\overline{f}$ and $(Tf)^{*}d$ are in $\mathcal{M}$ and so is their 
composite $\overline{d}$. The resulting partial map 
$\ptlmap{f^{*}y = (\overline{d},\overline{y})}{TX}{X}$ is a partial 
algebra, the domain of the lifting.

 As for universality, given any morphism of partial algebras 
$\morph{g}{x}{y}$
over $\monomorph{f}{X}{Y}$, we  consider the corresponding 
factorisation
throught the limit diagram:

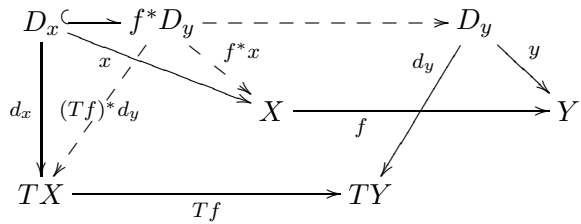
\[\xymatrixrowsep{1.5pc}
 \xymatrixcolsep{1.5pc}
\begin{diagram}  
{D_{x}}\ddto_-{d_{x}}\ar@^{(->}[r]\drrto_(.3){x}&\ar@{-->}[ddl]|-{(Tf)^{*}d_{y}}
{f^{*}D_{y}}\ar@{-->}[dr]^-{f^{*}x}\ar@{-->}[rrr]&&& 
\ddlto_(.3){d_{y}}{D_{y}}\drto^{y}&\\
 &&{X}\ar@{->}[rrr]_(.3){f} 
&&&{Y}\\
 {TX}\ar@{->}[rrr]_-{Tf}&&& {TY}&&   
\end{diagram}
\]
\noindent which stablishes the required cartesian property. 
As for saturation, it is clear that the isomorphism required
transfers via limits (pullbacks) from the codomain partial algebra to 
the domain one.
\qed
\end{proof}

A morphism of partial algebras $\morph{f}{x}{y}$ is an
\textbf{epimorphism} whenever
$\morph{f}{X}{Y}$ is so (in $\cat{B}$). Recall from \cite{Borceux94I} 
the notion of factorisation system.
 As a consequence of Proposition \ref{prop:fib-ptl-alg} we have the 
following:

\begin{proposition}\label{prop:factorisation-system}
    Any morphism $\morph{f}{x}{y}$ of relational algebras factors 
(uniquely) as an 
    epimorphism followed by a (monic) Kleene morphism. The pair of 
classes of morphisms
    (epimorphisms/monic Kleene
    morphisms) constitutes a pullback stable factorisation system in
    $\Ptl\mbox{\scriptsize -}\alg{T}$.
\end{proposition}
\begin{proof}
    Given a morphism of partial algebras $\morph{f}{x}{y}$, consider 
    its underlying morphism $\morph{f}{X}{Y}$ in $\cat{B}$ and take 
its 
    epi/$\mathcal{M}$-mono factorisation $f = m{\circ}e$ (the image 
factorisation of 
    \sectref{subsec.internal-envelope}). 
    Consider the cartesian lifting $\morph{\overline{m}}{m^{*}y}{y}$ 
    and the corresponding induced map of partial algebras over $e$,
    $\morph{e'}{x}{m^{*}y}$ (which follows by cartesianness, or 
simply 
    considering the limit diagram defining $m^{*}y$). This yields the 
    required factorisation $f = \overline{m}{\circ}e'$.
    Both (epimorphisms) and
    (Kleene morphisms) are clearly stable under pullback.\qed
\end{proof}

\begin{remark}
The reader familiar with factorisation systems in $\Cat$ would 
recognise the 
above factorisation as formally analogous to the 
bijective-on-objects/fully-faithful factorisation of functors.
\end{remark}

\begin{remark}
    The bicategory of partial maps is a subbicategory of 
    that of relations: the relational composition of two partial maps 
is 
    again one such.
    Partial algebras are thus a special instance of \textbf{relational
    algebras}: a relational algebra for a monad
    $\morph{T}{\cat{B}}{\cat{B}}$ (which preserves jointly-monic 
spans) is given by a relation $\rel{x}{TX}{X}$
    and 2-cell $\morph{\alpha}{x{\circ}Tx}{x{\circ}\mu}$ satisfying
    $x{\circ}(\mathit{id},\eta) = (\mathit{id},\mathit{id})$. 
Morphisms
    are defined as per partial algebras, yielding the category 
    $\Rel\mbox{\scriptsize -}\alg{T}$ of relational algebras. 
    The category of partial algebras is thus a full subcategory of 
$\Rel\mbox{\scriptsize -}\alg{T}$.

    The forgetful functor 
    $\morph{\underline{U}}{\Rel\mbox{\scriptsize -}\alg{T}}{\cat{B}}$
    is a fibration, whose cartesian liftings are precisely those
    morphisms such that the tent diagram is a limit. So the
    fibrational properties of $\Ptl\mbox{\scriptsize -}\alg{T}$ of
    Proposition \ref{prop:fib-ptl-alg} are
    inherited from this ambient fibred category.
\end{remark}

The construction in the proof of 
Corollary \ref{cor:reflectivity-saturated-ptl-alg} is an instance of 
the cartesian lifting of $\mathcal{M}$-monos to partial algebras of 
Proposition \ref{prop:fib-ptl-alg}. To see this, we must show that 
the 
unit $\morph{q\circ\eta_{X}}{X}{{\mathsf{E}}(x)}$ is a 
$\mathcal{M}$-monic. 
\begin{lemma}
    \label{lemma:unit-mono}
    If $\cell{\eta}{\mathit{id}}{T}$ is a 
\textit{cartesian\/}\footnote{The 
property of being cartesian implies that the unit is monic. It 
usually holds for free-monoid constructions.} 
trasnformation, that is, the corresponding naturality squares are 
pullbacks (or cartesian squares), whose instances are
$\mathcal{M}$-monos, the unit of the adjunction of Theorem 
\ref{cor:envelope-partial-algebras}
is an $\mathcal{M}$-mono.
\end{lemma}
\begin{proof}

We would show that the composite 
$q{\circ}\eta_{X}$ is the $\mathcal{M}$-image of $\eta_{X}$ along $q$.

Recall the definition of ${\mathsf{E}}(x)$ via the coequaliser

\[\xymatrixrowsep{1.5pc} \xymatrixcolsep{3pc}
    \let\objectstyle=\scriptstyle
    \let\labelstyle=\scriptscriptstyle
        \begin{diagram} 
    {TD}\rrcompositemap<-2>_{Td_{x}}^{\mu}{\omit} 
    \rrto^-{Tx}&&{TX}\ar@{->>}[r]^-{q}&{{\mathsf{E}}(x)}
           \end{diagram}
           \]
\noindent To compute the $\mathcal{M}$-image of (the 
$\mathcal{M}$-mono) $\eta_{X}$, we pull-back both morphisms defining 
the coequaliser along $\eta_{X}$, take their coequaliser, and obtain 
the image as the uniquely induced morphism between the coequalisers 
(\cf. the proof of Theorem   \ref{cor:envelope-partial-algebras}). 

We show that pulling back the pair $(Tx,\mu{\circ}Td)$ 
along $\eta_{X}$ yields a pair of equal morphisms, and thus its 
coequaliser is an isomorphism, and the corresponding image is the 
composite $q{\circ}\eta_{X}$ as desired. 
By cartesiannes of $\eta$, 
we have a pullback diagram
\[\xymatrixrowsep{1.5pc} \xymatrixcolsep{3pc}
    \let\objectstyle=\scriptstyle
    \let\labelstyle=\scriptstyle
        \begin{diagram} 
            {D}\dto_-{\eta_{D}}\rto^-{x}&{X}\dto^-{\eta_{X}}\\
    {TD}\rto_-{Tx}&{TX}
           \end{diagram}
           \]

\noindent For the other morphism, we have the pullback

\[\xymatrixrowsep{1.5pc} \xymatrixcolsep{3pc}
    \let\objectstyle=\scriptstyle
    \let\labelstyle=\scriptstyle
        \begin{diagram} 
          {X}\dto_{\eta'_{x}}\rto^{\mathit{id}}&  
{X}\dto_-{\eta_{X}}\rto^-{\mathit{id}}&{X}\ddto^-{\eta_{X}}\\
   {D}\dto_-{\eta_{D}}\rto^-{d}& {TX}\dto_-{\eta_{TX}}& \\
    {TD}\rto_-{Td}&{T^{2}X}\rto_-{\mu}&{TX}
           \end{diagram}
           \]
           
\noindent where the top-left square is a pullback by the axiom for 
the partial algebra $x$ (see \sectref{subsec.partial-alg-monad}), the 
bottom-left square is a pullback by 
cartesiannes of $\eta$, and the right rectangle is a pullback 
because $\eta_{X}$ is mono, and thus in the following diagram
\[\xymatrixrowsep{1.5pc} \xymatrixcolsep{3pc}
    \let\objectstyle=\scriptstyle
    \let\labelstyle=\scriptstyle
    \begin{diagram} 
      {X}\dto_{\mathit{id}}\rto^{\mathit{id}}&  
{X}\dto_-{\eta_{X}}\rto^-{\mathit{id}}&{X}\ddto^-{\eta_{X}}\\
   {X}\dto_-{\eta_{X}}\rto^-{\eta_{X}}& {TX}\dto_-{\eta_{TX}}& \\
    {TX}\rto_-{T\eta_{X}}&{T^{2}X}\rto_-{\mu}&{TX}
       \end{diagram}
       \]
the outer and top-left squares are pullbacks and  the bottom-left 
square 
is a pullback by cartesianness of $\eta$. Since $\mu$ is (split) epi 
and we are in a regular setting, we conclude the right rectangle is a
pullback indeed.

Thus we are left to coequalise the pair 
$\morph{(x{\circ}\eta'_{x},\mathit{id})}{X}{X}$ to obtain the image, 
but 
$x{\circ}\eta'=\mathit{id}$, again by the axiom for the partial 
algebra $x$.
\qed            
\end{proof}

\begin{remark}
In $\Set$, the composite $q\circ\eta_{X}$ being mono can be shown 
quite easily reasoning with elements:
$$(\tuple{z}/_{\sim} = \tuple{y}/_{\sim}) \Longrightarrow 
\tensor_{1}\tuple{z} = \tensor_{1}\tuple{y}) \Longrightarrow (z = y) 
$$
           \noindent and we could have deduced the above result from 
           this \textit{logical argument\/} by appealing to the full
       exact embedding
of a regular category in a Grothendieck topos \cf.\cite[\S 
2.7]{Borceux94II}.           
\end{remark}

F{i}nally, let us relate Kleene morphisms 
with  saturation of partial algebras. We have the composite
adjunction

\[\xymatrixrowsep{1.5pc} \xymatrixcolsep{3pc}
    \begin{diagram} 
    {\Ptl\mbox{\scriptsize 
-}\alg{T}}\ar@/^2ex/[rr]^-{{\mathsf{E}}}_-{\perp}&&\ar@/^2ex/[ll]^-{U}{\alg{T}_{P}}\ar@/^2ex/[rr]^-{{\underline{U}}}_-{\perp}&&\ar@/^2ex/[ll]^-{\mathsf{Tot}}{\alg{T}}
       \end{diagram}
       \]
\noindent where $\underline{U}$ forgets the distinguished subobject
of the algebra and its right adjoint takes the algebra
$\morph{x}{TX}{X}$ to the pair $(x,\monomorph{\mathit{id}}{X}{X})$
(the \textit{trivial\/} or \textit{total\/} distinguished subobject).
Let $\underline{T} = U{\mathsf{Tot}}\underline{U}{\mathsf{E}}$ the
resulting monad on $\Ptl\mbox{\scriptsize -}\alg{T}$ with unit
$\cell{\rho}{\mathit{id}}{\underline{T}}$.

\begin{corollary}\label{cor:factorisation-saturation}\hfill{ }\\
    Under the hypothesis of Lemma \ref{lemma:unit-mono},
    for a partial algebra $T$-algebra $\ptlmap{x}{TX}{X}$ tfae:
    \begin{enumerate}
        \item \label{uno} $\ptlmap{x}{TX}{X}$ is saturated

        \item\label{zwei}  $\morph{\rho_{x}}{x}{\underline{T}(x)}$ is 
a Kleene
    morphism
    
    \item\label{tres} $x$ embeds into a $T$-algebra via a Kleene 
morphism
    
\end{enumerate}
\end{corollary}
\begin{proof}
 \hfill{ }\\   
\noindent\underline{(\ref{uno}) $\Longrightarrow$ (\ref{dos})}:
We only need to point out that the partial algebra induced by a total
one $\morph{x}{TX}{X}$ and a given subob{j}ect 
$\monomorph{m}{D}{X}$ is the domain of the cartesian lifting of
$\mathsf{Tot}(x)$ along the $\mathcal{M}$-mono $m$. The argument in
the proof of Corollary \ref{cor:reflectivity-saturated-ptl-alg} shows
that when $x$ is saturated the unit $\tilde{\eta}_{x}\cong
\Sigma_{q}(d)$ is cartesian (Lemma \ref{lemma:unit-mono} shows it is 
in $\mathcal{M}$).

\noindent\underline{(\ref{zwei}) $\Longrightarrow$ (\ref{tres})}:
immediate

\noindent\underline{(\ref{tres}) $\Longrightarrow$ (\ref{uno})}:
Given a monic Kleene morphism $\morph{m}{x}{y}$ into a total algebra 
$y$, since $m$ is cartesian $x$ is saturated (Proposition
\ref{prop:fib-ptl-alg}).
\qed
\end{proof}

\begin{remark}
    The latter equivalence between saturation and $\rho_{x}$ being a
    Kleene morphism is  Freyd's purported statement  that 
    \textit{the envelope constructruction embeds every paracategory
    into a category via a Kleene functor\/} (the set-theoretic proof 
of this
    admittedly non-evident statement is unfinished in the
   draft manuscript \cite{Freyd96}). However, as we show in  
    \cite{HermidaMateus02b}, the (composite) adjunction between
    partial and ordinary algebras is not a useful tool to study the
    structure of the former; we must keep track of the distinguished
    subobjects to obtain meaningful results.    
\end{remark}

An interesting instance of the notion of partial algebra appears in 
our follow-up article \cite{HermidaMateus02b}:
\textit{partial multicategories\/}. For a 
partial multicategory we  formulate  
\textit{representability\/}, in the sense of \cite{Hermida99a} (for a 
restricted class of 
morphisms). A representable partial multicategory provides a suitable 
axiomatisation for certain 
`partial tensor products' which arise in the context of probabilistic 
automata \cite{Mateus2000}.

\appendix
\section{Background material}
\label{sec:preliminaries}

In this section we review the categorical background we use as
framework for our treatment of internal paracategories. F{i}rst, we
recall the relevant notions of monoidal category and their morphisms 
(\sectref{sec:monoidal}).
Secondly, monoidal categories give the ambient structure in which to 
def{i}ne monoids. Internal monoids in a monoidal category are 
classif{i}ed
by (strong) monoidal functors from the simplicial category 
(\sectref{sec:monoids}).
Specialising to the setting of endospans in a category with
pullbacks, we get the notion internal
category as an instance of that of internal monoid 
(\sectref{sec:internal-categories}). 

The context in which we apply (a lax version of) this general theory 
of internal monoids
is that of a bicategory (\sectref{subsec:bicat-lax-functors}) (or 
locally 
ordered category) of partial maps
(\sectref{sec:partial-maps}).

\subsection{Monoidal categories and their morphisms}
\label{sec:monoidal}

Our basic framework is that of a {\em monoidal category}, \ie. a 
category 
$\cat{V}$ endowed with a pseudo-associative 
$\morph{\tensor}{\cat{V}\times\cat{V}}{\cat{V}}$ and pseudo-unitary 
$I\in\cat{V}$ structure. See \cite{MacLane98,Kelly82} for details. 
The primary examples of interest to us are:

\begin{enumerate}
\item \label{One} $\Set$, the category of sets and functions, with 
monoidal structure given by cartesian products, 

\item \label{Two} $\Ptl_{\cal M}(\cat{B})$, the category of partial 
maps relative to a class of {\em monic domains\/} $\cal M$ in a 
category with f{i}nite limits $\cat{B}$ (which we recall in 
\sectref{sec:partial-maps} below), with monoidal structure given by 
f{i}nite products, and

\item \label{Three} $\Spn{\cat{B}}(C,C)$, the category of endospans 
on 
an object $C$ in a category $\cat{B}$ with pullbacks, 
\cf.\sectref{sec:internal-categories}.

\end{enumerate}

We recall the relevant notions for morphisms of monoidal categories:

\begin{itemize}
    \item  Given monoidal categories 
$(\cat{M},\tensor,I)$ and
$(\cat{N},\tensor',I')$, a \textbf{monoidal functor} 
between them is a functor $\morph{F}{\cat{M}}{\cat{N}}$ together with 
{\em structural natural transformations} 
\[ \morph{\gamma}{I'}{FI}\qquad 
\morph{\delta_{x,y}}{Fx{\tensor}'Fy}{F(x{\tensor}y)} \] 
\noindent for $x,y$ in $\cat{M}$ subject to 
coherence axioms, which guarantee that we get well-def{i}ned 
comparison 
transformations 
$\morph{\delta_{\vec{x}}}{\bigotimes_{1}^{n}Fx_{i}}{F(\bigotimes_{1}^{n}x_{i})'}$

for any  $n$-ary tensor of $\vec{x} = \tuple{x_{1},\ldots,x_{n}}$.

    \item  A monoidal functor is called \textbf{normal} when 
    $\morph{\gamma = \mathit{id}}{I'}{FI}$. It is called 
\textbf{strong} (resp. \textbf{strict}) when the structural natural 
transformations are isomorphisms (resp. identities).

    \item  Given monoidal functors $(F,\gamma,\delta)$ and 
    $(F',\gamma',\delta')$ a \textbf{monoidal transformation} between 
them 
    is a natural transformation $\cell{\alpha}{F}{F'}$ compatible 
    with the structural transformations, \ie.
    \[ \alpha_{I}\gamma = \gamma' \qquad 
    \alpha_{x{\tensor}y}\delta_{x,y} = 
    \delta'_{x,y}(\alpha_{x}\tensor\alpha_{y})\]

\end{itemize}

\begin{remark}
    A good way to understand laxity for monoidal categories is via 
{\em multicategories}, which yield an effective and simple 
description of classif{i}ers for lax functors, 
\cf.\cite[Rmk.~9.5]{Hermida99a}. The monoid classif{i}er $\Delta$ is
obtained as a special case of this construction.
\end{remark}

\subsection{Bicategories and lax functors}
\label{subsec:bicat-lax-functors}

The notions of monoidal functor and monoidal transformation have 
their several-objects 
    counterparts in the context of bicategories, where they become
\textbf{lax functors} and \textbf{lax transformations} respectively.

The notion of bicategory is part of the standard categorical tool-kit
\cite[\S XII.6]{MacLane98}. A bicategory $\cal K$ has objects, 
morphisms and 2-cells, so
that the morphisms ${\cal K}(X,Y)$ together with their 2-cells form a
category, and there are identities $\mathit{id}_{X}\in{\cal K}(X,Y)$
and composition functors $\morph{\tensor}{{\cal K}(X,Y)\times{\cal
K}(Y,Z)}{{\cal K}(X,Z)}$ together with unit
$\cell{\rho_{f}}{f}{\tensor(f,\mathit{id}_{X})}$,
$\cell{\lambda_{f}}{\tensor(\mathit{id}_{Y},f)}{f}$ (for
$\morph{f}{X}{Y}$) and associativity
$\cell{\alpha_{f,g,h}}{\tensor(\tensor(f,g),h)}{\tensor(f,\tensor(g,h))}$
 
2-cell isomorphisms (for composable $f,g,h$) satisfying coherence 
conditions as
per a monoidal category.

A \textbf{lax functor} $\morph{F}{{\cal K}}{{\cal L}}$
maps objects to objects, morphisms to morphisms and 2-cells to
2-cells, preserving the source-target relationships. But it is not
required to preserve the composition and identities. Instead, we have
\textit{comparison} 2-cells 
$\cell{\gamma_X}{\mathit{id}_FX}{F\mathit{id}_X}$ and 
$\cell{\delta_{f,g}}{Ff{\comp}Fg}{F(f{\comp}g)}$ satisfying the
evident coherence axioms (like those for monoidal functors).
When the comparison 2-cells are isomorphisms, we refer to the lax
functor as a \textbf{homomorphism of bicategories}.

A \textbf{lax transformation} 
    $\morph{\cell{\alpha}{F}{F'}}{{\cal K}}{{\cal L}}$ between lax
    functors assigns to every morphism 
    $\morph{h}{X}{Y}$ in $\cal K$ a 2-cell as displayed below
    
    \[\xymatrixrowsep{1.5pc}
 \xymatrixcolsep{1.5pc}
          \begin{diagram}   
{FX}\dlowertwocell<0>_{Fh}{^<-5>{{\alpha_h}}^{}}
\rto^-{\alpha_X}&{F'X}
\ar[d]^{F'h}\\
{FY}\ar[r]_-{\alpha_Y}&{F'Y}
\end{diagram}
\]
\noindent subject to coherence conditions which make them 
compatible with the comparison 2-cells of $F$ and $F'$.
More information on bicategorical matters (as relevant for this 
paper) appears in \cite{Hermida99a}, which provides further
references for the interested reader.

The two bicategories of interest in this paper are: that of partial
maps (\sectref{sec:partial-maps} below) and that of spans. Given a 
category $\cat{B}$
with pullbacks, we build the \textbf{bicategory of spans} 
$\Spn{\cat{B}}$: its objects are those of 
$\cat{B}$, morphisms are given by spans 
$X\stackrel{d_{R}}{\leftarrow}R\stackrel{c_{R}}{\rightarrow}Y$ and 
2-cells $\cell{f}{R}{S}$ correspond to 
arrows between the top objects of the spans such that the following 
diagram commutes:

  \[\xymatrixrowsep{1pc}
  \xymatrixcolsep{0.7pc}
        \begin{diagram}
          &\dlto_{d_R}{R}\ddto^{f}
          \drto^{c_R}& \\
          {X}& &{Y} \\
          &\ulto^{d_S}{S}\urto_{c_S}&
               \end{diagram}
        \]
        
        The identity morphism on $X$ is the span 
$X\stackrel{\mathit{id}}{\leftarrow}X\stackrel{\mathit{id}}{\rightarrow}X$ 

and composition is given by pullback:

$$\prooftree
{\xymatrixrowsep{2pc}
 \xymatrixcolsep{1.5pc}
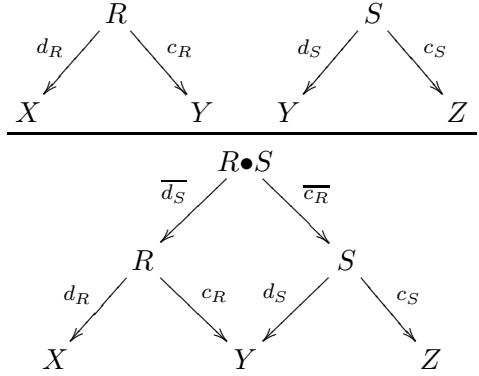
        \begin{diagram}
        &\dlto_{d_R}{R}\drto^{c_R}& &   &\dlto_{d_S}{S}\drto^{c_S}& \\
 {X}& &{Y}                         &  {Y}& &{Z} 
               \end{diagram}
}
\justifies
{\xymatrixrowsep{2pc}
 \xymatrixcolsep{1.5pc}
\begin{diagram}
&&\dlto_{\overline{d_S}}{R{\bullet}S}\drto^{\overline{c_R}} &&\\
&\dlto_{d_R}R\drto^{c_R}& &\dlto_{d_S}S\drto^{c_S}& \\
X& &Y& &Z \\
\end{diagram}
}
\thickness=0.08em
\endprooftree $$

\subsection{Internal monoids and their classif{i}er}
\label{sec:monoids}

A monoidal category provides the ambient structure to def{i}ne an 
{\em 
internal monoid\/}, so that for instance an internal monoid in $\Set$ 
is a monoid in the usual sense, while an internal monoid in 
\sectref{sec:monoidal}(\ref{Three}) above amounts to an {\em internal 
category}.

We recall from \cite[Ch.~VI]{MacLane98} the basic notion of monoid in 
a monoidal category and its classif{i}er.
\begin{definition}
  Given a monoidal category 
$\tuple{\cat{C},\tensor,I,\alpha,\rho,\lambda}$, a \textbf{monoid} in 
it consists of an ({\em underlying\/}) object $M$ and morphisms 
$\morph{e}{I}{M}$ ({\em unit\/}) and $\morph{m}{M{\tensor}M}{M}$ 
({\em multiplication\/}) satisfying the equations
  \begin{eqnarray*}
    m\comp(e\tensor\mathit{id}_M)\comp\lambda_M & = &  \mathit{id}_M 
\\
m\comp(\mathit{id}_M{\tensor}e)\comp\rho_M & = &  \mathit{id}_M \\
m\comp(m{\tensor}\mathit{id}_M) & = & 
m\comp(\mathit{id}_M{\tensor}m)\comp\alpha_{M,M,M}
  \end{eqnarray*}
A \textbf{monoid morphism} is a morphism between the underlying 
objects which commutes with unit and multiplication. We thus have the 
category $\Mon{\cat{C}}$ of monoids in $\cat{C}$.
\end{definition}

Given monoidal categories 
$\tuple{\cat{C},\tensor,I}$ and 
$\tuple{\cat{C}',\tensor',I'}$, a monoidal 
functor  $\morph{F}{\cat{C}}{\cat{C}'}$ (with 
structural cells $\gamma$ and $\delta$) induces a functor 
$\morph{\Mon{F}}{\Mon{\cat{C}}}{\Mon{\cat{C}'}}$ as follows:
 \[\xymatrixrowsep{1pc}
  \xymatrixcolsep{1pc}
        \begin{diagram}
          I\rrto^-{e}&&{M}\ar@{|->}[d]&&\llto_-{m}{M{\tensor}M} \\
  {I'}\rto_-{\gamma}& 
{FI}\rto_-{Fe}&{FM}&\lto^-{Fm}{F(M{\tensor}M)}&\lto^-{\delta_{M,M}}{FM{\tensor}FM}
\end{diagram}
        \]
Furthermore, notice that a monoid in $\cat{C}$ amounts to a monoidal
functor $\morph{M}{\bold{1}}{\cat{C}}$ and a monoid morphism 
to a  monoidal transformation between the corresponding monoidal 
functors. Thus the action of $\Mon{F}$ above amounts to composition 
with $F$. For a monoidal category $\cat{C}$ there is a universal 
construction of a monoidal category $G(\cat{C})$ such that monoidal  
functors out of $\cat{C}$ correspond to strong monoidal functors out 
of 
$G(\cat{C})$. Here we proceed to recall the explicit description of 
$G(\bold{1})$, which classif{i}es monoids \cf.\cite[\S 
VI.5]{MacLane98}.

\begin{definition}[Simplicial Category]
  Let $\Delta$ denote the category of f{i}nite ordinals 
$$[n] = \{ 0,\ldots,n-1 \} $$
(including the 
empty one $\emptyset = [0]$) and monotone functions between them. It 
carries a strict 
monoidal structure given by ordinal addition 
$\morph{+}{\Delta\times\Delta}{\Delta}$, with $\emptyset$ as unit. 
Since $[1] = \bold{1}$ is the terminal object, there are unique 
arrows 
$\morph{e}{\emptyset}{\bold{1}}$ and 
$\morph{m}{\bold{1}+\bold{1}}{\bold{1}}$, which satisfy (trivially) 
the equations for a monoid in $\Delta$.
\end{definition}

\begin{proposition}[Classif{i}cation of monoids]
 Given a monoidal category $\cat{C}$, the functor 
$\morph{\_(\bold{1})}{\MonCat(\Delta,\cat{C})}{\Mon{\cat{C}}}$ which 
evaluates a strong functor (resp. transformation) at the monoid 
$\bold{1}$, induces an equivalence of categories:
 \begin{displaymath}
   \MonCat(\Delta,\cat{C})\simeq\Mon{\cat{C}}
 \end{displaymath}
  
\end{proposition}

The above proposition means that all the objects in $\Delta$ are 
f{i}nite tensor powers of $\bold{1}$, \ie.
$[n] = \bigotimes^{n}\bold{1}$ and that all (non-identity) maps 
are generated from $e$ and $m$ by $+$ and 
composition. Thus $\Delta$ gives us a neat way to organise all the 
$n$-ary operations present in a monoid, as well as their equational 
theory.

\subsection{Internal categories}
\label{sec:internal-categories}

Given a category $\cat{B}$ with f{i}nite limits, we can def{i}ne an 
{\em 
internal category\/} or a {\em category object\/} in it. The most 
suitable way to phrase such def{i}nition for our purposes is via the
associated bicategory of spans (\sectref{subsec:bicat-lax-functors}).

\begin{definition}\label{def:internal-category}
    An \textbf{internal category} $\mathsf{C}$ in $\cat{B}$ is a 
monad in 
    $\Spn{\cat{B}}$, \ie. an object $C_{0}$ in $\cat{B}$ and a span 
$C_{0}\stackrel{d}{\leftarrow}C_{1}\stackrel{c}{\rightarrow}C_{0}$ 
    endowed with a monoid structure in $\Spn{\cat{B}}(C_{0},C_{0})$, 
given by 
 2-cells $\cell{\iota}{C_{0}}{C_{1}}$ (\textit{identities\/})
    and $\cell{m}{C_{1}{\comp}C_{1}}{C_{1}}$ (\textit{composition\/}) 
    satisfying associativity and unit laws.
    
    An \textit{internal functor\/} from $\mathsf{C}$ to $\mathsf{D}$ 
    consists of a pair of morphisms
    ($\morph{f_{0}}{C_{0}}{D_{0}}$, $\morph{f_{1}}{C_{1}}{D_{1}}$),
    commuting with domain, codomain, identity and multiplication
    maps.
    
    With the evident pointwise composition of internal functors we
   obtain the category $\Cat(\cat{B})$ of internal categories and
    internal functors in $\cat{B}$.
\end{definition}

A monoid in 
$\Spn{\cat{B}}(C_{0},C_{0})$ corresponds to a strong monoidal functor 
$\morph{C}{\Delta}{\Spn{\cat{B}}(C_{0},C_{0})}$.

\begin{remark}\label{rem:simplicial-object}
    An internal category in $\cat{B}$ amounts also to a 
f{i}nite-limit 
    preserving functor 
    $\morph{C}{(\Delta_{+})^{\mathit{op}}}{\cat{B}}$, where 
$\Delta_{+}$ 
    is the category of non-empty f{i}nite ordinals and monotone 
functions. 
    This combinatorial/topological point of view of a category, due 
to Grothendieck, is usually referred to as its 
    {\em nerve\/}. It is a (special kind of) {\em simplicial 
object\/} 
    in $\cat{B}$. But this alternative formulation would be more 
diff{i}cult to generalise to the context of partial maps to obtain 
paracategories, as we do here.
\end{remark}

\subsection{Monoidal bicategory of partial maps}
\label{sec:partial-maps}

 %

Given a category $\cat{B}$, we consider a class of monos $\cal M$ in 
it such that:

\begin{itemize}
    \item  isomorphisms are in $\cal M$

    \item  for $\monomorph{m}{P}{X}$ in $\cal M$ and 
$\morph{h}{Y}{X}$ 
    in $\cat{B}$, the pullback
    \[ \xymatrixrowsep{2pc}
 \xymatrixcolsep{1.5pc}
\begin{diagram}
    f^{*}P\rto^{\bar{f}}\ar@^{(->}[d]_-{f^{*}m}&P\ar@^{(->}[d]^-{m}\\
Y\rto_{f} &X 
\end{diagram} \]
exists and $\monomorph{f^{*}m}{f^{*}P}{Y}$ belongs to $\cal M$

    \item  $\cal M$ is closed under composition.
\end{itemize}

Such a class $\cal M$ is called a \textit{dominion \/} in
\cite{RobinsonRosolini88,Rosolini86}, which are standard references
for this subject. 
We call $\cal M$-{\bf subob{j}ect} an isomorphism class of monos in 
$\cal 
M$ (over their common codomain). We abuse notation and write 
$\monomorph{m}{P}{X}$ 
for the equivalence class of $m$ qua $\cal M$-subob{j}ect.

\begin{examples}
    \begin{itemize}
    \item  Consider $\cat{B}= \omega$-\textit{CPO\/}
 of partial orders with suprema of countable chains and monotone maps
 preserving them. A suitable class $\mathcal{M}$ of monos is given by
 the \textit{admissible\/} subobjects, which are those closed under
 suprema in the ambient cpo (these are used as they admit Scott's
 fixpoint induction principle). This is the traditional basic setting
 of domain theory.
 
        \item  Consider $\cat{B}=\Set^{\cat{C}^{\mathit{op}}}$ any
    presheaf topos. Any topology on it determines suitable
    classes of $\mathcal{M}$, \eg. the closed subobjects.
    \end{itemize}
\end{examples}

\begin{definition}\label{def:partial-maps}
    The bicategory of \textbf{partial maps} $\Ptl_{{\cal 
M}}(\cat{B})$ 
    consists of 
    
    \begin{description}
        \item[objects]  those of $\cat{B}$
    
        \item[morphisms]  a morphism $\ptlmap{(m,f)}{X}{Y}$ is given 
        by a span, so that $\monomorph{m}{P}{X}$ is an $\cal 
        M$-subob{j}ect and $\morph{f}{P}{Y}$ is a morphism in 
     $\cat{B}$ (the total part of the partial map). We also write $\ptlmap{h}{X}{Y}$, in which case
     its total part is the morphism $\morph{h}{P_{h}}{Y}$, with 
     $\monomorph{d_{h}}{P_{h}}{X}$ the corresponding $\cal 
        M$-subob{j}ect.
     
        \item[2-cells]  Given $\ptlmap{(m,f),(n,g)}{X}{Y}$, a 2-cell 
        between them is a morphism of the corresponding spans. Since 
there 
        is at most one such (because we are considering $\cal 
        M$-subob{j}ects), we have in fact a partial order on 
        $\Ptl_{{\cal M}}(\cat{B})(X,Y)$
    \end{description}
    
    Composition and identities are inherited from $\Spn{\cat{B}}$; 
the conditions on $\cal M$ ensure that the composition of two partial 
    maps qua spans yields a partial map.
\end{definition}

Notice that we have an embedding of $\cat{B}$ into $\Ptl_{{\cal 
M}}(\cat{B})$: 
\[ \morph{f}{X}{Y} \mapsto \ptlmap{(\mathit{id},f)}{X}{Y} \]

\noindent This embedding enjoys a universal property which 
characterises $\Ptl_{{\cal 
M}}(\cat{B})$ up to equivalence \cite{Hermida02a}.
In general, a partial map $\ptlmap{(m,f)}{X}{Y}$ is in the image of 
this embedding iff $\monomorph{m}{P}{X}$ is an isomorphism. In this 
case, $\ptlmap{(m,f)}{X}{Y}$ is called a {\em total map\/}. Notice 
that total maps are {\em maximal\/} with respect to the partial order 
in their hom-sets. 

For a more vernacular notation using elements, the expression
$$ f(x_1,\ldots,x_n) \succeq g(x_1,\ldots,x_n) $$
\noindent means that $\ptlmap{f\leq 
g}{M_1\times\ldots{\times}M_n}{M}$, while 
$$ f(x_1,\ldots,x_n) = g(x_1,\ldots,x_n) $$
\noindent means (Kleene) equality of partial maps (equal domains and 
values). With these notation, 
composition of morphisms amounts to substitution of terms for 
variables. We also use the \textit{is def{i}ned\/} $(\_)\downarrow$
predicate, so that $f(e)\downarrow$ means that $e$ is in the domain
of def{i}nition of the partial map $f$. 

\begin{remark}
 \label{Kleene-equality}   
Notice that the partial order in the hom-sets $\Ptl_{{\cal 
M}}(\cat{B})(X,Y)$ is the $\succeq$ relation (inclusion of domains
and equality of results over the smaller one), while the induced
equality is Kleene equality.
\end{remark}

We have the following cancellation property of total maps:

\begin{lemma}\label{lem:cancellation-total-map}
    Given partial maps $\ptlmap{(m,f)}{X}{Y}$ and 
    $\ptlmap{(n,g)}{Y}{Z}$, if the composite $(n,g)\comp(m,f)$ is 
total, then 
    $(m,f)$ is total as well.
\end{lemma}
\begin{proof}
    Consider the composite of spans $(n,g)\comp(m,f)$:
    \[ \xymatrixrowsep{2pc}
 \xymatrixcolsep{1.5pc}
\begin{diagram}
&&\ar@^{(->}[dl]_-{f^{*}n}f^{*}Q{\comp}Q\drto^{\overline{f}} &&\\
&\ar@^{(->}[dl]_-{m}P\drto^{f}& &\ar@^{(->}[dl]_-{n}Q\drto^{g}& \\
X& &Y& &Z \\
\end{diagram}
 \]
By hypothesis, the composite $m{\comp}f^{*}n$ is an isomorphism. 
Hence $m$ is a split epi, and thus an isomorphism.

\qed
\end{proof}

The stability conditions for the monos in $\cal M$ ensure that the 
f{i}nite product structure on $\cat{B}$ extends to a monoidal 
structure 
on $\Ptl_{{\cal M}}(\cat{B})$, by applying the product functor to the 
spans corresponding to the partial maps:

$$\prooftree
{\xymatrixrowsep{2pc}
 \xymatrixcolsep{1.5pc}
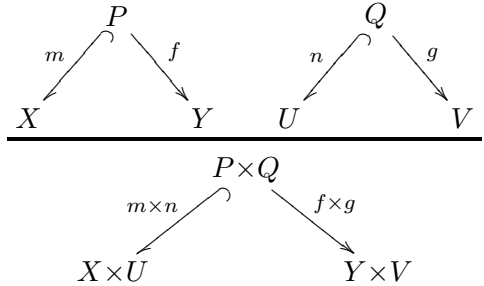
    \begin{diagram}
    &\ar@^{(->}[dl]_-{m}{P}\drto^{f}& &  
    &\ar@^{(->}[dl]_-{n}{Q}\drto^{g}& \\
 {X}& &{Y}                         &  {U}& &{V} 
           \end{diagram}
}
\justifies
{\xymatrixrowsep{2pc}
 \xymatrixcolsep{1.5pc}
\begin{diagram}
&\ar@^{(->}[dl]_-{m{\times}n}{P{\times}Q}\drto^{f{\times}g}&  \\
{X{\times}U}& &{Y{\times}V} \\
\end{diagram}
}
\thickness=0.08em
\endprooftree $$   

\noindent Furthermore, this monoidal 
structure is compatible with the partial order on the hom-sets. Thus 
$\Ptl_{{\cal M}}(\cat{B})$ is a {\em monoidal bicategory}, that is a 
bicategory (actually a 2-category) endowed with a monoid structure 
with respect to the cartesian product of bicategories\footnote{There 
is a more general notion of monoidal bicategory in the literature, 
where the monoidal structure is understood with respect to the 
Gray-tensor product of bicategories.}.

\comment{
\begin{remark}
\label{rmk:general-partial-maps}
    The def{i}nition of a category of partial maps above is standard. 
We note two immediate meaningful generalisations where partial 
algebraic structures should be explored:
    \begin{itemize}
        \item  The maps in the class $\cal M$ are assumed to be monos 
        as they represent the {\em domain of def{i}nition\/} of a 
        partial map. Of course, only the stability properties of 
        $\cal M$ are used in the def{i}nition of the bicategory 
        $\Ptl_{\cal M}(\cat{B})$ (which is forced into a category by 
        considering subob{j}ects). A suitable class $\cal M$ 
        with such properties would be the {\em formal\/} monos of a 
        stable factorisation system for $\cat{B}$. For instance in 
        $\Cat$ (the category of small categories) the {\em discrete 
        cof{i}brations\/} are the formal monos for a stable 
        factorisation system (whose formal epis are the so-called 
        initial functors). An interesting variant of this situation 
is to consider the anti-discrete cof{i}brations as monos: the 
resulting 
partial maps are Makkai's 
        {\em anafunctors\/} \cite{Makkai96}. Another case which has 
found its way into the literature is that where the class ${\cal M}$ 
consists of the {\em weak equivalences\/} of a model structure; the 
corresponding partial algebras (for operads) have been considered in 
\cite{KrizMay95}.
    
        \item  A partial map from $X$ to $Y$ is formally given by a 
        {\em property\/} on $X$ and a map from its {\em extent\/} to 
        $Y$. In constructive logic we would consider not mere 
        entailments between properties but suitable equivalence 
        classes of {\em proofs\/} between them and would therefore 
        have genuine categories of partial maps from $X$ to $Y$, 
        rather than mere partial orders. The general setting 
        to study this structure is a {\em f{i}bration with 
comprehension and 
        direct images\/}, the objects and arrows in the total 
category 
        corresponding to the `abstract' properties and their proofs 
of entailments respectively. 
    \end{itemize}
\end{remark}
}

\vskip5mm
\noindent{\underline{\textbf{Acknowledgements:}}

The authors gratefully acknowledge the prompt and comprehensive 
reviews by the 
referees. In particular, the second referee indicated a gap in our 
previous 
def{i}nition of internal paramonoid which led us to introduce the 
crucial saturation 
condition.

{\small
\bibliographystyle{alpha}
\bibliography{references}

}

\end{document}